\font\capit=cmcsc10 \font\addressit=cmcsc8 \font\eightrm=cmr8

\def\idfy{\mathbin{\hbox{$\widetilde\rightarrow$}}} \def\ext{{\rm Ext}}
\def\mult{{\rm mult}}  
\def\kum{{\rm Kum}} \def\ver{{\rm Ver}} 
\def\bl{{\rm Bl}} \def\cliff{{\rm Cliff}} \def\secant{{\rm Sec}}
\def\sectil{{\widetilde \secant}} \def\mult{{\rm mult}}

\def\eeetil{{\widetilde \varepsilon}} 
\def\dual{^{\vee}} 

    \def\g00{\bp \Gamma_{00}}
\def\trig{g^1_3} \def\trigg{h^1_3} \def\thsing{\Theta_{\rm sing}}

 \def\nm{{\rm Nm}} \def\pic{{\rm Pic}}

 \def\ff{{\cal F}}  \def\aa{{\cal A}}
\def\nn{{\cal N}}   \def\ll{{\cal L}}
 \def\tt{{\cal T}}  \def\ll{{\cal L}}
  \def\oo{{\cal O}} 
\def\ww{{\cal W}}

\def\su{{\cal SU}_C} \def\iic{{\cal I}_C} \def\eee{\varepsilon}

  \def\ses#1#2#3{0\rightarrow {#1}\rightarrow {#2}\rightarrow
{#3}\rightarrow 0}  
\def\z{{\bf Z}}   \def\c{{\bf C}}
\def\bp{{\bf P}}

\def\ad{{\rm ad}\ }

\def\im{{\rm im \ }}  \def\ker{{\rm ker \ }}
\def\rk{{\rm rank \ }}
\def\map#1{\ \smash{\mathop{\longrightarrow}\limits^{#1}}\ }

\def\pf{{\it Proof.}\ }

\def\<{\langle} \def\>{\rangle}

\newcommand{\qed}{{\unskip\nobreak\hfill\hbox{ $\Box$}\medskip\par}}

\documentstyle[epsf,leqno,bezier,12pt]{article}

\makeatletter                                             \renewcommand{\@begintheorem}[2]{                        
\sl \trivlist \item [\hskip \labelsep {\bf #2\ \ #1.}]}                        
\makeatother                                           
\makeatletter \def\section{\@startsection {section}{1}{\z@}{-3.5ex plus
-1ex minus
 -.2ex}{1.5ex plus .2ex}{\large\bf}}

\def\subsection{\@startsection{subsection}{2}{\z@}{-3.25ex plus -1ex
minus
 -.2ex}{1.5ex plus .2ex}{\normalsize\it}}

\newcommand{\numberequationsassubsubsections}

\newtheorem{prop}{Proposition}[section] \newtheorem{lemm}[prop]{Lemma}
\newtheorem{theo}[prop]{Theorem} \newtheorem{cor}[prop]{Corollary}
\newtheorem{rem}[prop]{\it Remark} \newtheorem{rems}[prop]{\it Remarks}
\newtheorem{ex}[prop]{Example}

\begin{document}

\title{Subvarieties of $\su(2)$ and $2\theta$-divisors in the Jacobian}
\author{W.M. Oxbury, C. Pauly and E. Previato}

\date{}

\maketitle

Let $\su(2,L)$ denote the projective moduli variety of semistable rank
2 vector bundles with determinant $L\in \pic(C)$ on a smooth curve $C$
of genus $g>2$; and suppose that $\deg L$ is even. It is well-known
that, on the one hand, the singular locus of $\su(2,L)$ is isomorphic
to the Kummer variety of the Jacobian; and on the other hand that when
$C$ is nonhyperelliptic $\su(2,\oo)$ has an injective morphism into the
linear series $|2\Theta|$ on the Jacobian $J_C^{g-1}$ which restricts
to the Kummer embedding $a \mapsto \Theta_a + \Theta_{-a}$ on the
singular locus. Dually $\su(2,K)$ injects into the linear series
$|\ll|$ on $J_C^0$, where $\ll = \oo(2\Theta_{\kappa})$ for any theta
characteristic $\kappa$, and again this map restricts to the Kummer map
$J_C^{g-1} \rightarrow |2\Theta|\dual = |\ll|$ on the singular locus.
This map to projective space (the two cases are of course isomorphic)
comes from the complete series on the ample generator of the Picard
group, and (at least for a generic curve) is an embedding of the moduli
space.  Moreover, its image contains much of the geometry studied in
connection with the Schottky problem; notably the configuration of
Prym-Kummer varieties.

In this paper we explore a little of the interplay, via this embedding,
between the geometry of vector bundles and the geometry of
$2\theta$-divisors. On the vector bundle side we are principally
concerned with the Brill-Noether loci $\ww^r \subset \su(2,K)$ defined
by the condition $h^0(E) >r$ on stable bundes $E$. These are analogous
to the very classical varieties $W_{g-1}^r \subset J_C^{g-1}$. Unlike
the line bundle theory, however, general results---connectedness,
dimension, smoothness and so on---are not known for the varieties
$\ww^r$ (see \cite{BF}).

On the $2\theta$ side we shall consider the Fay trisecants of the
$2\theta$-embedded Kummer variety, and the subseries $\g00 \subset
|\ll|$ consisting of divisors having multiplicity $\geq 4$ at the
origin. This subseries is known to be important in the study of
principally polarised abelian varieties \cite{vGvdG}: in the Jacobian
of a curve its base locus is the surface $C-C \subset J_C^0$ (plus a
pair of isolated points in the case $g=4$) \cite{W}, whereas for a ppav
which is not a Jacobian it is conjectured that the origin is the only
base point (but see \cite{BD}).

The organisation and main results of the paper are as follows. In the
first two sections we study two families of lines on $\su(2,K) \subset
|\ll|$ (or equivalently $\su(2,\oo) \subset |2\Theta|$), each of
dimension $3g-2$. These are the Hecke lines, coming from vector bundles
of odd degree, on the one hand, and lines lying inside $g$-dimensional
linearly embedded extension spaces (generating the lowest stratum of
the Segre stratification), on the other.  We prove ({\bf 1.3} and {\bf
1.4}) that a line in $\su(2,K)$ lies in both of these families if and
only if it intersects the Kummer variety. Moreover, we show ({\bf 2.1})
that every trisecant of the Kummer is such a line, and in particular
lies on the moduli space. (This fact is certainly well-known to the
experts, but we were not aware of a reference in the literature.) As a
corollary we show ({\bf 2.2}) that the Kummer variety has quadrisecant
lines if and only if the curve is hyperelliptic.

In section 3 we introduce the subschemes $\ww^r \subset \su(2,K)$, and
as a first step in their study we examine the stratification by $h^0$
of spaces of extensions, which will then map rationally into
$\su(2,K)$. The natural formulation of this stratification turns out to
involve the Clifford index, and as an easy by-product we obtain the
inequality ({\bf 3.7}) $$ h^0(E) \leq g+1 -\cliff(C) $$ for any
semistable rank 2 bundle $E$ with $\det E = K$.

In section 4 we prove, using a spectral curve construction, that ({\bf
4.1}) $$ \ww^2 = \g00 \cap \su(2,K) $$ provided $C$ is nonhyperelliptic
of genus 3 or 4, or nontrigonal of genus 5. On the other hand, we show
later on ({\bf 8.3}) that the equality fails for all curves of
genus~6.

The remaining four sections of the paper are devoted to examining some
of the geometry in detail for each of the cases $g= 3,4,5,6$.  For
genus~3 the moduli space $\su(2,K)$ is embedded in $\bp^7$ as the
unique Heisenberg-invariant quartic singular along the Kummer
variety---the so-called {\it Coble quartic}. We examine the
configuration $\ww^1 \subset \ww \subset \bp^6$, where $\bp^6$ is the
hyperplane spanned by the `generalised theta divisor' $\ww=\ww^0$. It
is known that $\ww$ has a unique triple point $\ww^2=\g00$; we show
({\bf 5.3}) that $\ww^1$ is a Veronese cone (already to be found in the
classical literature \cite{C}) with vertex $\ww^2$, and whose
generators are trisecants of the Kummer variety corresponding to a
natural embedding of $|K|$ in the parameter space of all trisecants. In
addition, we identify the tangent cone of $\ww$ at the triple point
({\bf 5.5}): this is nothing but the secant variety of the Veronese
surface, with equation $\det A =0$ where $A$ is a symmetric $3\times 3$
matrix.

To each nonhyperelliptic curve of genus 4 one can associate a nodal
cubic threefold $\tt\subset \bp^4$, which can be described in various
ways. The view we adopt here is that it is the rational image of
$\bp^3$ via the linear system of cubics containing the canonical curve.
There is an identification $\bp^4 \idfy \g00$, due to Izadi \cite{I},
and we prove ({\bf 6.4}) that this restricts to an isomorphism $\tt
\idfy \ww^2$, with the node mapping to $\ww^3 = \{\trig \oplus
\trigg\}$, the direct sum of the two trigonal line bundles on the
curve. For the proof of this we make use of Izadi's description of the
lines in $\tt$ as pencils of $2\theta$-divisors. Note also that by
passing to the tangent cone at the origin, $\g00$ may be viewed as a
linear system of quartics in canonical space $\bp^3$. We observe ({\bf
6.11}) that projection of the cubic $\ww^2$ away the node can be
naturally identified with the quotient of $\Gamma_{00}$ by $q^2$ where
$q$ is the unique quadric vanishing on the canonical curve.

For genus 5 we show ({\bf 7.2}) that $\ww^3$ is a Veronese surface
cutting the Kummer variety in the image of a plane quintic. If the
curve is nontrigonal this quintic is the discriminant of the net of
quadrics containing the canonical curve, while in the trigonal case it
is isomorphic to the projection of the canonical curve from a
trisecant.

Finally we show ({\bf 8.1}) that for a nontrigonal curve of genus 6 the
locus $\ww^4$ is a single point, stable if $C$ is not a plane quintic
(this case was also observed by Mukai in \cite{Muk}), while $\ww^4$ is
a line not meeting the Kummer in the trigonal case. In the generic case
$\ww^4$ is the vertex of a configuration of five $\bp^6$s which form
the intersection of $\su(2,K)$ with $\g00$ residual to $\ww^2$.

\medskip

{\it Acknowledgements.} The authors are grateful to R. Donagi, E. Izadi, P.
Newstead, S. Ramanan, E. Sernesi, and especially B. van Geemen for
various helpful comments. We would also like to thank Miles Reid for
his organisation of the Warwick Algebraic Geometry Symposium 1995--96
where much of this work was carried out, and the MRC for its
hospitality during this symposium; also L. Brambila-Paz and the
hospitality of IMATE-Morelia, Mexico during the Vector Bundles workshop
in July 1996.  Research of the third named author was partly supported
by NSA grant MDA904-95-H-1031.

\vfill\eject \bigskip\noindent {\large\bf I Lines}

\section{The $g$-plane ruling}

For a line bundle $L$ on the curve $C$ we denote by $\su(2,L)$ the
projective moduli variety of semistable rank 2 vector bundles with
determinant $L$; and in particular we shall be concerned with
$\su(2,K)$. The semistable boundary of this space is the image of
$J_C^{g-1} \rightarrow \su(2,K)$ mapping $L\mapsto L\oplus KL^{-1}$;
this is the singular locus when $g>2$. Throughout the paper we shall
view both the moduli space and the Kummer variety---when $C$ is
nonhyperelliptic---as lying in the projective space $|2\Theta|\dual =
\bp^{2^g-1}$ in the standard way: by the complete linear series $|\ll|$
on $\su(2,K)$, where $\ll$ is the ample generator of $\pic\ \su(2,K)
\cong \z$, restricting to $\oo(2\Theta)$ on the Jacobian. A {\it line}
on $\su(2,K)$ is then an embedded $\bp^1$ on which the restriction of
$\ll$ has degree one.

We shall consider the following subvariety of $\su(2,K)$ ruled by
$g$-planes. For $x\in \pic^{g-2}(C)$ let $\bp(x) = \bp H^1(C,K^{-1}x^2)
\cong \bp^g$. This parametrises isomorphism classes of extensions $$
\ses{x}{E}{Kx^{-1}}, $$ and thus has a moduli map to $\su(2,K)$, which
is linear (with respect to $\ll$) and injective (see also remark
\ref{bertram} below). Globally we have a ruling:  $$ \begin{array}{rcl}
\bp U & \map{\eee}& \su(2,K) \\ &&\\
 \downarrow \bp^g&&\\ &&\\ \pic^{g-2}(C) &&\\ \end{array} $$ where $U =
R^1\pi_* K^{-1}\nn^{2}$, with $\nn \rightarrow C\times \pic^{g-2} (C)$
a Poincar\'e bundle and $\pi: C\times \pic^{g-2}(C) \rightarrow
\pic^{g-2}(C)$ the natural projection.

We shall need to make repeated use, in what follows, of the following
result of Lange--Narasimhan \cite{LN}. Consider any extension $$
\ses{n_0}{F}{n_0^{-1}\otimes \det F} $$ where $n_0 \subset F$ is a line
subbundle of maximal degree. This is represented by a point $f$ of the
extension space $\bp H^1(C, n_0^2 \otimes \det F\dual) = \bp H^0(C,
Kn_0^{-2}\otimes \det F)\dual$, into which the curve $C$ maps via the
linear series $|Kn_0^{-2}\otimes \det F|$. For an effective divisor $D$
on $C$, we shall denote by $\overline D$ the linear span in $\bp H^1(C,
n_0^2 \otimes \det F\dual)$ of the image of this divisor. Then the
following is proposition 2.4 of \cite{LN}:

\begin{lemm} \label{ln} With the above notation there is a bijection,
given by $\oo (D) = n^{-1}n_0^{-1}\otimes \det F$, between:
\begin{enumerate} \item line subbundles $n\subset F$, $n\not= n_0$, of
maximal degree; and \item line bundles $\oo (D)$ with degree $\deg D =
\deg F - 2 \deg n_0 $ and such that $f\in \overline D$.
\end{enumerate} \end{lemm}

Let us return now to the $g$-planes $\bp(x)\hookrightarrow \su(2,K)$,
where $x\in \pic^{g-2}(C)$, and the following well-known facts.  The
curve $C$ maps into $\bp(x)$ via $|K^2x^{-2}|$, as a special case of
the Lange--Narasimhan picture.  Moreover, a point of $\bp(x)$
represents a stable bundle (with $x$ as maximal line subbundle)
precisely away from the image of $C$; while a point $q\in C\subset \bp
(x)$ represents the equivalence class of the semistable bundle
$x(q)\oplus Kx^{-1}(-q)$. In other words there is a commutative
diagram:  \begin{equation} \label{5} \begin{array}{rcl} C & \map{t_x} &
J_C^{g-1} \\ &&\\ |K^2x^{-2}|\downarrow && \downarrow\\ &&\\ \bp(x)
&\map{\eee}& \su(2,K)\\ \end{array} \end{equation} where $t_x:q\mapsto
x(q)$ and the second vertical arrow maps $L\mapsto L\oplus KL^{-1}$.

The incidence relations between $g$-planes of this ruling can be given
as follows.

\begin{prop} \label{40} Suppose that $C$ is nonhyperelliptic.  For $x,y
\in \pic^{g-2}(C)$ the intersection $\bp(x) \cap \bp(y)$ is either
empty, or:  \begin{enumerate} \item the secant line $\overline{pq}$ of
the curve (in either of $\bp(x)$ or $\bp(y)$) if $x\otimes y =
K(-p-q)$; \item the point $x(p) \oplus Kx^{-1}(-p) \in
\kum(J_C^{g-1})$, if $H^0(C, Kx^{-1} y^{-1}) =0$ and $x\otimes y^{-1} =
\oo(q-p)$.  \end{enumerate} \end{prop}

\pf First we note that at any point $E\in \bp (x)$ away from the curve
the residual $g$-planes $\bp(y)$ through $E$ can be identified, via
lemma \ref{ln}, with the set of effective divisors $p+q$ such that $E$
lies on the secant line $\overline{pq}$; and the line bundles $x,y$ are
then related by \begin{equation} \label{30} x\otimes y = K(-p-q).
\end{equation} Note that the point $p$ on the curve in $\bp(x)$
represents the bundle $$ x(p)\oplus Kx^{-1}(-p) = y(q)\oplus
Ky^{-1}(-q), $$ i.e. it coincides with the image of $q$ on the curve in
$\bp(y)$; and similarly $q\in C\subset \bp(x)$ coincides with $p\in
C\subset \bp(y)$. This shows that condition (\ref{30}) is equivalent to
$\overline{pq} \subset \bp(x) \cap \bp(y)$.

On the other hand, when $C$ is nonhyperelliptic $\bp(x)$ and $\bp(y)$
cannot intersect in dimension greater than one: for then a generic
point $E$ of the intersection would lie on distinct secant lines
$\overline{pq}$ and $\overline{rs}$, both satisfying (\ref{30}), and
hence $\oo(p+q) = \oo(r+s)$, a contradiction.

The only other possibility for nonempty intersection $\bp(x) \cap
\bp(y)$ is that this intersection is a single point, in which case it
must be a point of the Kummer, and we easily find case 2.  \qed

Next we recall the Hecke correspondence between $\su(2,K)$ and
$\su(2,K(p))$ where $p\in C$ is a point of the curve. For a stable
bundle $F\in \su(2,K(p))$ we shall write $l_F \cong \bp^1 \subset
\su(2,K)$ for the image of $$ \begin{array}{rcl} \bp {\rm Hom}(F,\c_p)
&\rightarrow& \su(2,K)\\
 \phi &\mapsto& \ker \phi.\\ \end{array} $$ This is called the {\it
Hecke line} associated to the bundle $F$.  Our aim in the remainder of
this section will be to compare Hecke lines in $\su(2,K)$ with lines
contained in the $g$-planes $\bp(x)$.

\begin{theo} \label{10} Let $x\in \pic^{g-2}(C)$.  A projective line
$l\subset \bp(x)$ is Hecke if and only if it meets the image of the
curve $C\subset \bp(x)$.  \end{theo}

\pf Suppose that a Hecke line $ l_F$ is contained in $\bp(x)$. This
means that for all extensions of sheaves of the form $$
\ses{E}{F}{\c_p}, $$ where $K(p)$ is the determinant of $F$, the kernel
$E$ contains $x$ as a line subbundle. This means we have a pencil of
homomorphisms $x \hookrightarrow F$, for which we have two
possibilities.

{\it Either} the image subsheaf is constant and is in the kernel of
every homomorphism $F\rightarrow \c_p$. Then there is an inclusion of
sheaves $x(p) \subset F$, which by stability of $F$ is a line
subbundle, i.e. $F$ is an extension $$ \ses{x(p)}{F}{Kx^{-1}}.  $$ But
the space $\bp H^1(C,K^{-1} x^{2}(p))$ of such extensions is the image
of the projection of $\bp (x) = \bp H^1(C, K^{-1} x^{2})$ from the
point $p\in C \subset \bp (x)$, i.e. the set of lines in $\bp (x)$
passing through the point $p$. It is not hard to see that the line $l$
corresponding to $F$ in this manner is precisely $l_F$---and we note
that {\it every} line meeting the curve arises in this way.

{\it Or}---the second possibility---the image sheaf is non-constant, in
which case we have a subsheaf $x\oplus x\subset F$, with torsion
quotient supported on some effective divisor $D$. But then $\det F =
K(p)$ implies that $\oo(D) = Kx^{-2} (p)$. So $\deg D = 3$ and we
observe that $$ h^0(C, K^2x^{-2}(-D)) = h^0(C,K(-p)) = g-1, $$ i.e.
that $\dim \overline{D} = 1$. In this case the Hecke line $l_F$ is just
$\overline{D}$ and is trisecant to the curve.  \qed

It follows by a dimension count that for $g>2$ a generic Hecke line in
$\su(2,K)$ is not contained in any $g$-plane. We shall show next that
those that are are precisely the Hecke lines that meet the Kummer
variety (i.e. the singular locus).

\begin{theo} \label{20} Let $l_F \subset \su(2,K)$ be any Hecke line.
\begin{enumerate} \item There is a canonical surjection from line
subbundles $n\subset F$ with $\deg n = g-1$ to points of intersection
$l_F \cap \kum(J_C^{g-1})$ with the Kummer variety, which is bijective
if $l_F$ is not a tangent line of the Kummer.  \item The intersection
$l_F \cap \kum(J_C^{g-1})$ is nonempty if and only if $l_F$ is
contained in a $g$-plane $\bp(x)$ for some $x\in \pic^{g-2}(C)$. If $C$
is nonhyperelliptic and $l_F \cap \kum(J_C^{g-1})$ has cardinality $k$
then the number of such $g$-planes is $1+{k\choose 2}$.
\end{enumerate} \end{theo}

\begin{rems}\rm \label{21} 
{\it (i)} We have two irreducible
families of lines in $\su(2,K)$: the Hecke lines and the lines
contained in $g$-planes of the ruling. These families have the same
dimension $3g-2$, and we expect that each is a component of the Hilbert
scheme of all lines. Theorems \ref{10} and \ref{20} would then say that
{\it the intersection of these two components consists of the members
of each family which meet the Kummer variety}.  

{\it (ii)} The
cardinality of $l_F \cap \kum(J_C^{g-1})$ satisfies $k\leq 4$ if $C$ is
hyperelliptic and $k\leq 3$ otherwise; this follows from part 1 and
proposition 5.1 of \cite{Lan}. See also corollary \ref{4-secants}
below.  \end{rems}

{\it Proof of theorem \ref{20}:} For $n\subset F$ consider the
diagram:  $$ \begin{array}{rrcl} &&0&\\ &&\downarrow&\\ &&Kn^{-1}&\\
&&\downarrow&\\ 0\rightarrow n \rightarrow F & \rightarrow&
Kn^{-1}(p)&\rightarrow 0\\ &\phi \searrow&\downarrow&\\ &&\c_p&\\
&&\downarrow&\\ &&0&\\ \end{array} $$ Then if $E=\ker \phi$ we have a
surjective sheaf map $E\rightarrow Kn^{-1}\rightarrow 0$, and hence
$n\subset E$. So $E$ is S-equivalent to $n\oplus Kn^{-1}$ and defines a
point of intersection $l_F \cap \kum(J_C^{g-1})$.  To see that this is
surjective let $n\oplus K n^{-1} \in \kum(J_C^{g-1})$ be a point of
intersection with the Hecke line $l_F$. This means there is an exact
sequence $$ \ses{E}{F}{\c_p} $$ where $E$ is S-equivalent to $n\oplus
Kn^{-1}$; i.e. at least one of $n$ or $Kn^{-1}$ is a line subbundle of
$E$, and hence of $F$.

In a moment we shall verify that in this construction we have:
\begin{equation} \label{tangent} n\oplus Kn^{-1} \subset F \quad
\Longleftrightarrow \quad \hbox{$l_F$ is tangent to $\kum(J_C^{g-1})$
at $n\oplus Kn^{-1}$.} \end{equation} This will show that the
correspondence is bijective when $l_F$ is not a tangent line.

For part 2, we note first that if $l_F$ is contained in a $g$-plane
$\bp(x)$ then by theorem \ref{10} it meets the curve $C\subset \bp(x)$
and hence the Kummer. For the converse it suffices, by the preceding
construction, to suppose that there is a degree $g-1$ line subbundle
$n\subset F$.  Then letting $x = n(-p)$, it follows that $F$ is
represented by a point of the extension space $\bp H^1(C, K^{-1}n^2
(-p)) = \bp H^1(C, K^{-1}x^{2} (p))$ and hence determines---as in the
proof of theorem \ref{10}---a line $l \subset \bp (x)$ through $p\in
C$, which coincides with the Hecke line $l_F$.

Notice that by lemma \ref{ln} any residual degree $g-1$ line subbundles
$m\subset F$ correspond to points $q\in C$ by the relation $K(p-q) =
n\otimes m$; and in this case $l_F$ must be the secant $\overline{pq}
\subset \bp(x)$. In particular, $Kn^{-1}$ is a subbundle of $F$ if and
only if $l_F$ is the tangent line to $C\subset \bp(x)$ at $p$---this
proves (\ref{tangent}).

Finally, if we fix any $g$-plane containing $l_F$ then by proposition
\ref{40} the residual such $g$-planes correspond bijectively to the
effective divisors $p+q$ such that $l_F = \overline{pq}$, i.e. to pairs
of intersection points of $l_F$ with the Kummer. And so we obtain
$1+{k\choose 2}$ for the number of such $g$-planes.  \qed

\section{Trisecants of the Kummer variety} \label{trisecants}

Recall that the quotient $\kum(J_C^{g-1})$ of $J_C^{g-1}$ by the Serre
involution is embedded in $\bp^{2^g -1}$ by the linear system
$|2\Theta|$, and that this embedding extends to the moduli space
$\su(2,K)$ containing the Kummer as its singular locus (when $g>2$). In
this embedding the Kummer possesses a unique irreducible 4-dimensional
family of trisecant lines, which characterise Jacobians amongst
principally polarised abelian varieties. We shall briefly describe this
family of trisecants (see~\cite{M} or~\cite{Deb}).

The base $\ff$ of the family is the fibre product:  $$
\begin{array}{rcl} \ff &\rightarrow & S^4 C\\ &&\\ \downarrow
&&\downarrow \hbox{Abel-Jacobi}\\ &&\\ \pic^{g-3}(C) & \rightarrow &
\pic^4(C)\\ \end{array} $$ where the bottom map sends $a\mapsto K
a^{-2}$.  An element of $\ff$, in other words, is a pair $(a,D) \in
\pic^{g-3}(C)\times S^4 C$ such that $a^{2} = \oo(K-D)$. Writing $D=
p+q+r+s$, one shows that the following three points of $\bp^{2^g -1}$
are collinear:  \begin{equation} \label{3points} \begin{array}{lllll}
\phi(a(q+r)) && \phi(a(p+r)) && \phi(a(p+q))\\ =\phi(a(p+s)), &&
=\phi(a(q+s)), &&=\phi(a(r+s)).\\ \end{array} \end{equation}

We shall refer to the lines of $\bp^{2^g -1}$ parametrised by $\ff$ in
this way as the {\it Fay trisecants}.

\begin{theo} \label{fay} The Fay trisecants are precisely the Hecke
lines which are trisecant to the Kummer variety.  In particular they
all lie on $\su(2,K)$.  \end{theo}

\pf We ask for the condition on a Hecke line $l_F \subset \su(2,K)$ for
it to be trisecant to the Kummer. Let $\det F = K(p)$. For $l_F$ to
meet the Kummer $F$ must have a line subbundle $n \subset F$ of degree
$g-1$. Then by theorem \ref{20}, $l_F$ is a trisecant if and only if
$F$ has two further degree $g-1$ line subbundles. By stability these
are maximal, and so by lemma \ref{ln} the residual subbundles
correspond bijectively to points of $C$ mapping to the extension class
of $F$ under $$ C\ \map{|K^2 n^{-2}(p)|}\ \bp H^1(C,K^{-1} n^2(-p)).
$$ Thus trisecants $l_F$ correspond to {\it nodes} of the image curve
under the linear series $|K^2 n^{-2}(p)|$; and the condition for such a
node is that for points $q,r \in C$, $$ \begin{array}{rcl} h^0(C, K^2
n^{-2}(p-q-r) & \geq & h^0(C, K^2 n^{-2}(p))-1\\ &=& g-1;\\ \end{array}
$$ or equivalently $h^0(C,K^{-1} n^2(-p+q+r)) \geq 1$. This in turn is
equivalent to $K^{-1}n^2(-p+q+r) = \oo(s)$ for some $s\in C$. We
conclude that the necessary and sufficient condition for $l_F$ to be a
trisecant $\overline{pqr}$ of $C \subset \bp(x)$, $x=n(-p)$, is:
\begin{equation} \label{34} n^2 = K(p-q-r+s) \quad \hbox{or
equivalently} \quad x^2 = K(-p-q-r +s).  \end{equation} One can now
verify, using (\ref{5}), that the points of intersection of $l_F$ with
the Kummer---i.e. with the curve $C\subset \bp (x)$---are the three
points (\ref{3points}) where $a= x(-s) $.  \qed

\begin{cor} \label{4-secants} \begin{enumerate} \item If $C$ is
nonhyperelliptic then no Fay trisecant has more than three intersection
points with the Kummer variety.  \item If $C$ is hyperelliptic then all
Fay lines are exactly quadrisecant.  \end{enumerate} \end{cor}

\pf By theorems \ref{20} and \ref{fay} every Fay line lies in $\bp(x)$
for some $x\in \pic^{g-2}(C)$. Let $D$ be an effective divisor on $C$
with $\dim \overline{D} = 1$ in $\bp(x)$. This is equivalent to
$h^0(C,K^2x^{-2}(-D)) = g-1$, or, by Riemann-Roch, $$
h^0(C,K^{-1}x^{2}(D)) = \deg D - 2.  $$ If $\deg D = 5$ then this says
that $|K^{-1}x^{2}(D)|$ maps $C$ birationally to a plane cubic, which
is impossible; while if $\deg D = 4$ then it is equivalent to:
\begin{equation} \label{35} K x^{-2} = \oo(D-H) \end{equation} where
$H$ is a hyperelliptic divisor. This proves part 1; for part 2 let
$\overline{pqr}\subset \bp(x)$ be the trisecant constructed in the
proof of theorem \ref{fay}, and consider $D = p+q+r+\tau(s)$ where
$\tau:C\leftrightarrow C$ is the hyperelliptic involution. Then
(\ref{35}) follows from (\ref{34}) and we see that $\overline{pqr} =
\overline{D}$ is a quadrisecant.  \qed

\vfill\eject \bigskip\noindent {\large\bf II Loci}

\section{Brill-Noether loci in $\su(2,K)$}

Let $\ww \subset \su(2,K)$ be the closure of the locus of stable
bundles $E$ for which $H^0(C,E) \not= 0$, i.e. the `theta divisor' for
rank 2 bundles. In terms of the map $\phi: \su(2,K) \rightarrow
|2\Theta|\dual$, $\ww$ is the unique hyperplane section tangent to the
Kummer variety $\kum(J^{g-1}_C)$ along the image of the theta divisor
$\kum(\Theta)$.

The Brill-Noether loci are the subschemes $$ \ww \supset \ww^1 \supset
\cdots \supset \ww^{g-1} \supset \ww^g\\ $$ where $\ww^r$ is the
closure of the set of stable bundles $E$ for which $h^0(C,E) \geq
r+1$.  (We shall see in a moment that $\ww^{g+1} = \emptyset$---see
proposition \ref{bound}).

\begin{rem}\rm \label{petri} The local structure of $\ww^r$ is governed
by a symmetric Petri map $$ S^2 H^0(C,E) \rightarrow H^0(C,K\otimes \ad
E), $$ where $\ad E$ is the bundle of trace-free endomorphisms; as a
consequence $\ww^r$ has expected codimension $r+2 \choose 2$ in
$\su(2,K)$. (See for example \cite{BF}.)

In addition, it is not hard to show that {\it $\ww^1$ is the union of
all Hecke lines meeting $\ww^2$}.  We shall see this illustrated for
curves of genus 3 in theorem~\ref{coblecone}.  \end{rem}

In order to study the Brill-Noether loci $\ww^r$ we shall analyse them
first in spaces $\bp \ext^1(K-D,D)$ of extensions \begin{equation}
\label{basicext} \ses{\oo(D)}{E}{\oo(K-D)} \end{equation}

\noindent for some line bundle $\oo(D)\in \pic^d(C)$.  Usually, though
not always, we shall think of $D$ as an effective divisor; indeed $E$
has sections if and only if it can be expressed as such an extension
with $D$ effective.

\begin{rems}\rm \label{nag} 
{\it (i)} Note that by semistability
$d \leq g-1$ with equality if and only if $E$ is S-equivalent to
$\oo(D) \oplus \oo(K-D)$. Moreover, every $E\in \su(2,K)$ is such an
extension for some $D\in \pic^d(C)$ with $$ d\geq \Bigl[ {g-1 \over 2}
\Bigr].  $$ This follows from a classical result of Segre and Nagata
(see \cite{Lan2}) which says that every ruled surface of genus $g$ has
a section with self-intersection at most $g$.  

{\it (ii)} It will
be convenient below to introduce the Clifford index $\cliff(D) = \deg D
- 2r(D)$, where $r(D)=h^0(D) -1$, into our notation. Recall that the
Clifford index $\cliff(C)$ of the curve is defined to be the minimum
value of $\cliff(D)$ for which $h^i(D)\geq 2$ for $i=0,1$ (see
\cite{GL}).  Recall also that $$ \cliff(C) \leq \Bigl[ {g-1 \over 2}
\Bigr] $$ with equality for generic $C$.  \end{rems}

As in section 1, the curve $C$ maps into the space of such extensions;
and we shall denote the rational coarse moduli map of this space by
$\eee_D$:  $$ C\ \map{|2K-2D|}\ \bp \ext^1(K-D,D)\cong \bp^{3g-4-2d}
\ \map{\eee_D}\ \su(2,K).  $$ We shall denote by $\iic$ the ideal sheaf
of the image curve in $\bp^{3g-4-2d}$; and we shall write $\bl_C$ for
the blow-up along this curve and $\secant^n C$ for the variety of its
$n$-secant $(n-1)$-planes; although of course the map $C \rightarrow
\bp \ext^1(K-D,D)$ is not necessarily an embedding or even birational.

We shall write $$ \ww_D = \eee_D (\bp \ext^1(K-D,D)).  $$ Note here
that by $\eee_D (\Omega)$, where $\Omega \subset \bp \ext^1(K-D,D)$, we
shall always mean the proper transform of $\Omega$, i.e. the closure in
$\su(2,K)$ of the image of the domain of definition of $\eee_D$.

\begin{rem}\rm \label{bertram} The rational map $\eee_D$ has been
studied in detail by Bertram and others (see \cite{B}), and resolves to
a morphism $\eeetil_D$ of the blow-up:  $$ \begin{array}{r} \bp
\ext^1(K-D,D) \leftarrow \bl_C \leftarrow \bl_{\sectil_2 C} \leftarrow
\cdots \leftarrow \bl_{\sectil_{g-2-d}C} \\ \\ \downarrow \eeetil_D\\
\\ \su(2,K)\\ \end{array} $$ Moreover, hyperplanes of $|\ll| =
|2\Theta|\dual$ pull back to divisors of $|\iic^{g-2-d}(g-1-d)|$ on
$\bp \ext^1(K-D,D)$. We shall need these facts only in the cases $d=
g-2$ (already discussed in section 1) and $d=g-3$; in both of these
cases $\eee_D$ comes from the complete series $|\iic^{g-2-d}(g-1-d)|$.
\end{rem}

It is easy to analyse the filtration of each $\bp \ext^1(K-D,D)$ by the
dimension $h^0(E)$. For any such extension we have

\begin{equation} \label{h^0(E)} \begin{array}{rcl} h^0(E) &=& h^0(D) +
h^0(K-D) - \rk \delta(E)\\ &=& g+1 - \cliff(D) - \rk \delta(E)\\
\end{array} \end{equation} where $\delta (E): H^0(K-D) \rightarrow
H^1(D)$ is the coboundary homomorphism in the cohomology sequence of
(\ref{basicext}). (Note that (\ref{h^0(E)}) gives an upper bound
$g+1-\cliff(C)$ on $h^0(E)$; see proposition \ref{bound} below.)

By Serre duality $\delta(E)$ is an element of $\otimes^2 H^1(D)$, while
its transpose $\delta(E)^t$ is the coboundary map for the dual sequence
tensored with $K$:  $ \ses{\oo(D)}{K\otimes E\dual}{\oo(K-D)}.  $ But
$K\otimes E\dual = E$ and so $\delta (E) = \delta(E)^t$. We have
therefore constructed a linear homomorphism \begin{equation} \delta :
\ext^1(K-D,D) \rightarrow S^2 H^1(D), \end{equation} with respect to
which $h^0(E)$ satisfies (\ref{h^0(E)}).  But the {\it rank}
stratification of $S^2 H^1(D)$ coincides with the {\it secant}
stratification of its embedded Veronese variety $$ \begin{array}{rcl}
\ver : \bp H &\hookrightarrow& \bp S^2 H\\
       \xi &\mapsto& \xi\otimes \xi,\\ \end{array} $$ where $H =
H^1(D)$. In other words $$ \secant^{n} (\ver \bp H) = \{\ a \in \bp S^2
H\  |\  \rk a \leq n\ \} $$ for $n= 1,\ldots ,\dim H = r(D) -d+g$.  On
the other hand, the homomorphism $\delta$ is dual to the multiplication
map $ S^2 H^0(K-D) \rightarrow H^0(2K-2D) $ and so the above Veronese
embedding fits into the following commutative diagram:

\begin{equation} \label{verdiagram} \begin{array}{rcl} C & \map{|K-D|}
& \bp H^1(D) \\ &&\\ \scriptstyle|2K-2D| \displaystyle\downarrow
&&\downarrow \ver \\ &&\\ \bp \ext^1(K-D,D) &\map{\delta} & \bp S^2
H^1(D) \\ &&\\ \eee_D \downarrow &&\\ &&\\ \ww_D &&\\ \end{array}
\end{equation} We now define:  $$ \begin{array}{rcl} \Omega^0_D &=& \bp
\ker \delta, \\ \Omega^n_D &=& \delta^{-1} (\secant^n(\ver \bp
H^1(D))), \quad n= 1,\ldots , g-d+r(D).\\ \end{array} $$ (When it is
convenient we shall drop the subscript and write $\Omega^n
=\Omega^n_D$.) Thus if $\Omega^0$ is nonempty then $\Omega^0 \subset
\Omega^1 \subset \cdots \subset \Omega^{g-d+r(D)} \subset \bp
\ext^1(K-D,D)$ is a sequence of {\it cones} with vertex $\Omega^0$.

We can therefore state the main conclusion of this section as follows:

\begin{equation} \label{conethm} h^0(E) = g+1-\cliff(D)-n \qquad
\hbox{for $E\in \Omega_D^n \backslash \Omega_D^{n-1}$}.  \end{equation}

\begin{ex} \label{d=g-2} $\bf d=g-2.$ \rm If $D\in \pic^{g-2}(C)$ then
$\bp (D) = \bp \ext^1(K-D,D)$ is a $g$-plane of the ruling of section
1. The map $\eee_D: \bp(D) \hookrightarrow \su(2,K)$ is a linear
embedding, and we shall not distinguish between $\bp (D)$ and its
image. In this case (\ref{conethm}) says:  $$ \ww_D \cap \ww^r =
\Omega^{2r(D) + 2 -r}\subset \bp(D), \qquad r= r(D), \ldots , 2r(D) +
2.  $$ The cones $\Omega^n$ are constructed using $\delta :
\ext^1(K-D,D) \rightarrow S^2 H^1(D)$ where $\dim \ext^1(K-D,D) = g+1 $
and $\dim H^1(D) = h^0(D) + 1$.

If $h^0(D) =0$ then $\dim S^2 H^1(D) = 1$ and $\Omega^0 \subset
\Omega^1 = \bp(D)$ is a hyperplane, on which $h^0(E) = 1$. In other
words $\Omega ^0 = \bp(D) \cap \ww$.

If $h^0(D) = 1$ then $\dim S^2 H^1(D) = 3$. Since $|K-D|$ is a pencil,
$S^2 H^0(K-D) \rightarrow H^0(2K-2D)$ is necessarily injective and so
again $\delta$ is surjective. In this case, therefore:  $$
\begin{array}{rclr} \Omega^0 &\cong& \bp^{g-3} &\hbox{on which $h^0(E)
= 3$,}\\ \Omega^1 &=& \hbox{quadric of rank 3,} &h^0(E) = 2,\\ \Omega^2
&= & \bp(D),& h^0(E) =1.\\ \end{array} $$

If $h^0(D) = 2$ then the series $|K-D|$ maps $f: C\rightarrow \bp^2$
with degree $g$; and the homomorphism $\delta$ is no longer surjective
in general.  In fact surjectivity fails precisely when $f$ maps $C$
onto a conic, which cannot happen if $g$ is odd, but can occur for a
trigonal curve of genus 6, for example: if $|L| = g^1_3$ then take
$D=K-2L$.  In case $\delta$ {\it is} surjective we have:  $$
\begin{array}{rclr} \Omega^0 &\cong& \bp^{g-6} &\hbox{on which $h^0(E)
= 5$,}\\ \Omega^1 &=& \hbox{cone over a Veronese}&\\ &&\hbox{surface in
$\bp^5$,} &h^0(E) = 4,\\ \Omega^2 &=& \hbox{cone over cubic} &\\
&&\hbox{hypersurface $S^2 \bp^2 \hookrightarrow \bp^5$,} &h^0(E) = 3,\\
\Omega^3 &= & \bp(D),& h^0(E) =2.\\ \end{array} $$ And so on.  \end{ex}

The `universal' case of (\ref{conethm}) is the case $D=0$. This says
that $\ww^r$ is composed of the image of $\Omega^{g-r} \subset \bp
\ext^1(K,\oo)\cong \bp^{3g-4}$ together with those of the corresponding
cones in the exceptional divisors of the blow-up of
remark~\ref{bertram}:  $$ \ww^r = \bigcup_{D\geq 0\atop \deg D \leq
g-2} \eee_{D} \Omega^{g-r-\cliff(D)}.  $$ Diagram (\ref{verdiagram})
becomes in this case:

\begin{equation} \label{verdiagram2} \begin{array}{rcccl} &&C&
\map{|K|} & \bp^{g-1}\\ &&&&\\ &&\scriptstyle|2K|
\displaystyle\downarrow  && \downarrow \ver\\ &&&&\\
\Omega^{g-r}&\subset & \bp^{3g-4}&\map{\delta}& \bp S^2 H^1(\oo)\\
&&&&\\ \downarrow &&\eee_0\downarrow &&\\ &&&&\\ \ww^r & \subset& \ww
&&\\ \end{array} \end{equation}

\begin{rem}\rm Note that if $C$ is nonhyperelliptic then by Noether's
theorem $\delta$ in (\ref{verdiagram2}) is injective. Then $\Omega^n$
is the intersection of $\bp^{3g-4}\subset \bp S^2 H^1(\oo)$ with the
secant variety $\secant^n (\ver \bp^{g-1})$, and in particular contains
$\secant^n C\subset \bp^{3g-4}$. One can show, in fact, that Green's
conjecture on the syzygies of the canonical curve (see \cite{GL})
implies:  $$ \Omega^n = \secant^nC \subset \bp^{3g-4} \qquad \hbox{for
$n<\cliff(C)$.} $$ One consequence of this statement is that
$\Omega^{\cliff(C)}$ is the smallest cone in the sequence containing
semistable extensions.  \end{rem}

We conclude this section with two inequalities.  The first, which will
be useful later, is due to Mukai (\cite{Muk}, proposition 3.1):

\begin{lemm} \label{muk} If $|D|$ is base-point-free then for any rank
2 vector bundle $E$ with $\det E = K$ we have $ h^0(E(-D)) \geq h^0(E)
- \deg D.  $ \end{lemm}

\begin{prop} \label{bound} For all semistable bundles $E$ in $\su(2,K)$
we have $$ h^0(E) \leq g+1-\cliff(C).  $$ \end{prop}

\begin{rem}\rm In particular, this bound becomes:  $$ h^0(E)\leq
\cases{g & for nonhyperelliptic $C$,\cr
       g-1 & for $C$ not trigonal or a plane quintic,\cr \cdots&\cr
      [g/2] + 2& for generic $C$.\cr} $$ The bound $h^0(E) \leq g$ for
nonhyperelliptic curves was observed by Laszlo \cite{L}, proposition
IV.2.  \end{rem}

{\it Proof of proposition \ref{bound}.} We may assume $E$ comes from an
extension in $\bp \ext^1(K-D,D)$ where, by remark \ref{nag} (i) $$
\Bigl[ {g-1 \over 2} \Bigr] \leq \deg D \leq g-1.  $$ The right-hand
inequality implies $h^0(D) \leq h^1(D)$ so that if $h^0(D)\geq 2$ then
$\cliff(D) \geq \cliff(C)$ by definition. If, on the other hand,
$h^0(D) \leq 1$ then $\cliff(D) \geq \deg D \geq \cliff(C)$ by the
left-hand inequality together with remark \ref{nag} (ii). In either
case, therefore, the proposition follows from (\ref{h^0(E)}).  \qed

\section{$\g00$}

We shall as usual identify $\su(2,K)$ with its image in $|2\Theta|\dual
= |\ll|$ where $\ll = \oo(2\Theta_{\kappa}) \in \pic(J_C)$ for any
theta characteristic $\kappa \in \vartheta(C)$. Namely, a stable bundle
$E\in \su(2,K)$ is mapped to the divisor $D_E \in |\ll|$ defined by $$
D_E = \{ \ L\in J_C\ |\ h^0(C,L\otimes E)>0\ \}.  $$ On the other hand,
one can consider the linear system $\g00\subset |\ll|$ defined by:

\begin{equation} \begin{array}{rcl} \g00 &=& \{\ D\in |\ll|\ |\ \mult_0
D \geq 4\ \}\\
     &=& \{\ D\in |\ll|\ |\ C-C \subset {\rm supp}\ D\ \}.\\
\end{array} \end{equation} For the equivalence of these two definitions
see \cite{vGvdG} or \cite{W}; one can show, in addition, that $\g00$
has codimension $1+{1\over 2}g(g+1)$.

It is easy to verify that the Brill-Noether locus $\ww^2$ is always
contained in the subspace $\g00$. The main result of this section is a
partial converse:

\begin{theo} \label{w2g00} $\ww^2 \subset  \g00 \cap \su(2,K) \subset
\ww^1$. Moreover, if $C$ is nonhyperelliptic of genus 4 or nontrigonal
of genus 5 then $\ww^2 = \g00 \cap \su(2,K)$.  \end{theo}

\begin{rem}\rm We shall show later that $\ww^2 \not=  \g00 \cap
\su(2,K)$ for curves of genus 6 (see remark \ref{21dimspan}). For genus
4 the embedding $\ww^2 \subset \g00$ will be described in theorems
\ref{w2cubic} and \ref{quartics}.

\end{rem}

\begin{lemm} \label{4.3} Suppose $F$ is a semistable vector bundle of
rank 2 and degree $2d$ where $0\leq d \leq g-1$; and $k\geq 0$ an
integer. Then $h^0(F) \geq k$ if and only if $h^0(F(D)) \geq k $ for
all $D\in S^{g-1-d} C$.  \end{lemm}

Before proving this lemma let us show how it implies theorem
\ref{w2g00}. We suppose that $C-C \subset D_E$ for a stable bundle
$E\in \su(2,K)$, and we show first that $h^0(E)\geq 2$: by hypothesis
$h^0(E(p-q)) \geq 1$ for all $p,q\in C$, so by the lemma we deduce that
$h^0(E(-q)) \geq 1$ for all $q\in C$. From this it follows that $h^0(E)
\geq 2$, since $h^0(E) \geq 1$ and equality would imply that every
section vanishes at arbitrary $q\in C$, a contradiction.

Now suppose that $h^0(E) = 2$ and consider the evaluation map $e_q :
H^0(E) \rightarrow E_q$ for $q\in C$. Since $h^0(E(-q)) \geq 1$ we have
$\rk e_q \leq 1$ for all $q\in C$, and hence the sections of $E$
generate a line subbundle $L\subset E$ with $h^0(L) =2$. But by
stability of $E$ this must satisfy $\deg L < g-1$ so $C$ admits a
$g^1_{g-2}$.  So if $C$ is nonhyperelliptic of genus 4 or is
nontrigonal of genus 5 we obtain a contradiction, and we conclude that
$h^0(E) \geq 3$.  \qed

For the lemma, we shall  prove the following equivalent statement. Let
$E\in \su(2,\oo)$ and $\xi \in \pic^d(C)$, $0\leq d \leq g-1$. Then

\begin{equation} \label{killingD} H^0(C,\xi \otimes E)\geq k \quad
\Longleftrightarrow \quad {H^0(C,\xi(D) \otimes E)\geq k \atop
\forall\,\, D\in S^{g-1-d}C.} \end{equation}

We shall introduce a {\it spectral curve} $q: B= B_s\rightarrow C$ (see
\cite{BNR}). Namely, $B_s$ is the subscheme of the total space of the
canonical line bundle $K \map{q} C$ with equation $x^2 = s$, where
$s\in H^0(C,K^2)$ is a generic section. This is a smooth double cover
of $C$ of genus $g_B = 4g-3$, and there is a dominant rational map of
finite degree of the Prym variety on to $\su(2)$:  $$
\begin{array}{rcl} Q_s = \nm_q ^{-1}(K) &\rightarrow& \su(2)\\ \zeta
&\mapsto & q_*\zeta.\\ \end{array} $$ Moreover, the images of these
rational maps cover the moduli space as the section $s$ varies, and so
for any $E\in \su(2)$ we may assume that $E=q_* \zeta$ for some line
bundle $\zeta\in Q_s$, for suitable $s\in H^0(C,K^2)$.

By the projection formula the left-hand side of (\ref{killingD}) is $$
H^0(C,\xi \otimes E) = H^0(B_s, L) \qquad \hbox{where $L = \zeta
\otimes q^*\xi$.} $$ Notice that for $d<g-1$, $\deg L = 2g-2+2d \leq
4g-6 = g_B -3$, and in particular the Serre dual linear series $|K_B
L^{-1}|$ is base-point-free for generic $\zeta\in Q_s$. By choosing $s$
generically we may assume, for any given $E\in \su(2)$ and $\xi \in
J^d$, that this is the case.

We shall need:

\begin{lemm} Suppose, for $q:B\rightarrow C$ a double cover as above,
that $|N|$ is a base-point-free linear series on $B$. Then {\it either}
$N=q^*N'$ for some $N'\in \pic\ C$ {\it or} $$ h^0(B, N\otimes
q^*\oo(-x)) = h^0(B,N) - 2 $$ for generic $x\in C$.  \end{lemm}

\pf Write $q^{-1}(x) = x_1 + x_2$ with $x_1 \not= x_2$. Then $h^0(B,N)
- h^0(B,N(-x_1-x_2)) $ is the rank of the evaluation map $H^0(B,N)
\rightarrow \c_{x_1} \oplus \c_{x_2}$; and either this rank is 2 for
generic $x\in C$ or it is $\leq 1$ for all $x\in C$. In the latter
case, the base-point-free hypothesis ensures that the image of the
evaluation map is not contained in either summand; this implies that
every divisor in $|N|$ is symmetric, so $N= q^*N'$ as asserted.  \qed

We shall want to apply the lemma to $N = K_B L^{-1}$; we begin by
observing that this line bundle cannot be symmetric, as follows. Since
$K_B = q^* K_C^2$ and $L= \zeta \otimes q^*\xi$, $N =q^*N'$ would imply
that $\zeta = q^* \tau$ for some $\tau \in \pic\ C$. But then $$ E= q_*
\zeta = \tau \otimes q_* \oo_B = \tau \oplus K_C^{-1} \tau, $$
violating semistability.

So finally, consider the right-hand side of (\ref{killingD}). By the
projection formula this space is $$ H^0(C,\xi(D) \otimes E) = H^0(B,
L\otimes q^*\oo(D)).  $$ We note that $\deg L\otimes q^*\oo(D) = 4g-4 =
g_B -1$, so by Riemann-Roch $$ h^0(C,\xi(D) \otimes E) = h^0(B, K_B
L^{-1}\otimes q^*\oo(-D)).  $$ We now apply the lemma $e=g-1-d$ times
to $N = K_B L^{-1}$: this gives, for $D\in S^e C$ generic, $
h^0(C,\xi(D) \otimes E) = h^0(B, K_B L^{-1}) - 2e$.

\medskip {\it Proof of (\ref{killingD}).} Assuming the right-hand side
we have, by the last remark and by choosing $D$ generically, $h^0(B,
K_B L^{-1}) \geq k + 2e$. Consequently $$ \begin{array}{rcl} h^0(C,\xi
\otimes E) &=& h^0(B,L)\\ &=& h^0(B, K_B L^{-1}) + \deg L - g_B + 1 \\
&\geq& k+ 2e +\deg L  - g_B + 1 \\ &=& k.\\ \end{array} $$ The converse
is trivial.  \qed

\vfill\eject \bigskip\noindent {\large\bf III Low genera}

\section{Genus 3}

In this section we shall take $C$ to be a nonhyperelliptic curve of
genus 3. Then $\su(2,K)$ is the {\it Coble quartic} associated to the
Kummer variety in $\bp^7$ (see~\cite{NR} and \cite{C}, \S33).  It is
well-known that in this case the 3-plane ruling $\eee : \bp U
\rightarrow \su(2,K)$ of section 1 is surjective and has degree~8.
This follows easily from remark \ref{nag} (i) and lemma \ref{ln}.

The behaviour of $h^0(E)$ in each 3-plane of the ruling is given by
example~\ref{d=g-2}. Namely, if $h^0(x)=0$ then $\ww\subset \su(2,K)$
cuts $\bp(x)$ transversally in a 2-plane along which $h^0(E) = 1$,
while:

\begin{equation} \label{g=3,w1} \ww^1 = \bigcup_{p\in C} \Omega^1_p,
\qquad \Omega^1_p = \hbox{quadric cone} \subset \bp(p).  \end{equation}
In a moment we shall show that the vertices $\Omega^0_p$ of these cones
all coincide at a single point of $\su(2,K)$ (see \ref{w2point} and
\ref{coblecone}).

\begin{rem} \label{x=etap} \rm Note that for each $p\in C$ the image of
$C$ in $\bp(p)$ lies on the cone $\Omega^1_p$; and projecting along the
generators is the trigonality $f:C\rightarrow \bp^1$ given by the
series $|K(-p)|$.  (Conversely, one may show that the image of $C$ in a
3-plane $\bp(x)$ lies on a quadric cone only if $x= \eta(p)$ for some
$p\in C$ and some square root of the trivial line bundle, $\eta^2 =
\oo$.) \end{rem}

We consider now the birational map $\eee_0 : \bp \ext^1(K, \oo) = \bp^5
\rightarrow \ww \subset \su(2,K)$ and the diagram (\ref{verdiagram2}).
The map $\delta$ is an isomorphism: its dual $S^2 H^0(K) \rightarrow
H^0(2K)$ is surjective by Noether's theorem and injective since the
canonical curve $C\subset \bp^2$ is not contained in any quadric. Thus
$\Omega^1\subset \bp^5$ is the Veronese surface; and it is well-known
that its variety of secant lines $\Omega^2$ is a cubic hypersurface
isomorphic to $S^2 \bp^2$. Thus we have a tower of rational maps, where
$\eee_0$ is given by the complete linear series $|\iic(2)|$ on $\bp^5$
(see remark~\ref{bertram}):

$$ \begin{array}{rcl} \bl_C(\bp^5)&&\\ \downarrow&\searrow&\\ \bp^5 =
\Omega^3 &\map{\eee_0} & \ww\\ |&&|\\ S^2 \bp^2 = \Omega^2
&\longrightarrow& \ww^1\\ |&&|\\ C\subset \ver(\bp^2)= \Omega^1
&\longrightarrow& \ww^2\\ \end{array} $$ First of all, this allows us
to recover the following result of Laszlo \cite{L} and
Paranjape--Ramanan \cite{PR}. The bundle $V$ appearing here is simply
the normal bundle of $C$ canonically embedded in its Jacobian.

\begin{prop} \label{w2point} $\ww^2$ consists of a single point, i.e.
there is a unique stable bundle $V\in \su(2,K)$ with $h^0(V) =3$.
\end{prop}

\pf Since the canonical curve $C\subset \bp^2$ has degree 4 any quadric
of the series $|\iic(2)|$ either contains the Veronese surface
$\Omega^1$ or has no further points of intersection.  Thus $\Omega^1$
contracts to a single point $V\in\su(2,K)$ under~$\eee_0$.

On the other hand, for each $p\in C$ the 3-plane $\bp(p) \subset
\su(2,K)$ is the image of the fibre $\bp N_{C/\bp^5}$ of the
exceptional divisor in the blow-up, by remark \ref{bertram}. Making
this identification the point $\Omega^0_p\in \bp(p)$ is the normal
direction of $\Omega^1=\ver(\bp^2) \supset C$, and is therefore
contained in the closure of the image of $\Omega^1$.  In other words
$\Omega_p^0 = \{V\}$, and since we've seen that there are no further
points of $\ww^2$, this completes the proof.  \qed

We now wish to give a geometric description of $\ww^1 \subset
\su(2,K)$, and to this end we consider again the space $\ff$ of Fay
trisecants of the Kummer.  Notice that for genus 3 there is a natural
inclusion of the canonical series $$ |K| \hookrightarrow \ff
\ \map{J[2]}\ S^4 C $$ given by $D\mapsto (0,D)$ (see section
\ref{trisecants}). For $D\in |K|$ let us denote the corresponding
trisecant by $t_D \subset \su(2,K)$.  If $D=p+q+r+s$ then by the proof
of theorem \ref{fay}, $t_D$ lies in the four 3-planes
$\bp(p),\bp(q),\bp(r),\bp(s)$; in $\bp(p)$, $t_D$ is the trisecant line
$\overline{qrs}$, and similarly in the other three spaces.  By remark
\ref{x=etap}, on the other hand, this line is a generator of the cone
$\Omega^1_p \subset \bp(p)$, and conversely every generator is such a
trisecant. By (\ref{g=3,w1}), therefore, we conclude that

\begin{equation} \ww^1 = \bigcup_{D\in |K|}t_D \subset \su(2,K).
\end{equation}
From this we obtain the following description of $\ww^1$.

\begin{theo} \label{coblecone} The subvariety $\ww^1\subset
\su(2,K)\subset \bp^7 $ has the following structure:  \begin{enumerate}
\item $\ww^1$ is a cone over the Veronese surface $|K| = \bp^2
\map{|\oo(2)|} \bp^5$; \item $\ww^1$ has point vertex $\ww^2\in
\su(2,K)$; \item $\ww^1$ intersects the Kummer variety in the theta
divisor $\kum(\Theta)$, and projection along the generators of the cone
coincides with the 3 to 1 Gauss map $\kum(\Theta) \rightarrow |K|$.
\end{enumerate} \end{theo}

\pf We have already seen that each trisecant $t_{p+q+r+s}$, where
$p+q+r+s \in |K|$, passes through the point $V\in \ww^2$.  Assigning to
the divisor $p+q+r+s$ the tangent direction of $t_{p+q+r+s}$ at $V \in
\su(2,K)$ therefore defines an injective map $$ \pi : |K| \rightarrow
\bp T \su(2,K)|_{ V} \cong \bp^5; $$ for which $\pi^*\oo(1) = \oo(2)$
on the pencils $|K(-p)|\subset |K|$, and hence on the whole plane.
Parts 1 and 2 of the theorem follow straightaway.

From (\ref{5}) we see that the trisecant $t_{p+q+r+s}$ meets the Kummer
in the three points
$$ \begin{array}{c} \oo(p+s)\oplus \oo(q+r), \\ \oo(q+s)\oplus
\oo(p+r), \\ \oo(r+s)\oplus \oo(p+q); \\ \end{array} $$ which is
equivalent to part 3.  \qed

\begin{rem}\rm The Veronese cone of theorem \ref{coblecone} appears in
the work of Coble (\cite{C}, \S48). In particular, Coble exhibits a
uniquely determined cubic hypersurface in $\bp^6$ which cuts out
$\kum(\Theta)$ on the cone $\ww^1$.  It would be interesting to
interpret this cubic in terms of vector bundles.  \end{rem}

Finally, we shall sketch two proofs of the following fact.

\begin{theo} At the triple point $\ww^2= \{V\}$ the theta divisor $\ww$
has projectivised tangent cone $ \bp T_V \ww \cong \Omega^2 = S^2 \bp^2
\subset \bp^5 $.  \end{theo}

{\it First proof.} Since $V$ is stable we can identify $T_V \su(2,K)$
with $H^1(C,\ad V)$. We have already remarked (\ref{petri}) that since
$\det V = K$ the Petri map factors through the symmetric product $S^2
H^0(C,V) \rightarrow H^0(C,K\otimes \ad V)$; and this is dual to a map
$$ \mu: H^1(C,\ad V) \rightarrow S^2 H^0(C,V)\dual \subset {\rm
Hom}(H^0(V),H^1(V)).  $$ By standard Brill-Noether type arguments the
tangent cone $T_V \ww$ is the pull-back under $\mu$ of the
homomorphisms with nontrivial kernel (see for example \cite{L}). On the
other hand, one can show that $\mu$ is an isomorphism in the present
case. The tangent cone is therefore precisely the locus of singular
quadratic forms on $H^0(C,V)$, and hence isomorphic to $\Omega^2 = S^2
\bp^2$.  \qed

{\it Second proof} (due to B. van Geemen).  This exploits the fact that
the hypersurface $\ww \subset \bp^6$ has degree 4 (since it is a
hyperplane section of the Coble quartic), while $V$ is a triple point
of $\ww$ (\cite{L} proposition IV.7). It follows that $\bp T_V \ww $ is
the complement in $\bp^5$ of the (Zariski open) image of $\ww$ under
the rational projection map away from the point $V$.

We consider, then, the following diagram:

\begin{equation} \begin{array}{crccl} &&\ww&\subset & \bp^6\\ &&&&\\
&\hidewidth\eee_0 = \lambda_{|\iic(2)|}\nearrow&&\searrow&\downarrow
\pi_{V}\\ &&&&\\ \bp^5&&\map{\Delta}&&\bp^5\\ \end{array}
\end{equation} We have seen that $\ww$ is the (closed) image of $\bp^5$
under the rational map $\eee_0$ given by the complete linear series of
quadrics through the bicanonical curve, contracting $\ver(\bp^2)$ down
to the point $V$. Thus the rational map $\Delta$ is given by the
complete linear series of quadrics containing $\ver(\bp^2)$. It is
well-known that this can be identified with the inversion map of
symmetric $3\times 3$ matrices (geometrically, it sends a plane conic
to its dual conic). $\Delta$ is a birational involution, blowing up the
locus $\Omega^1 = \ver(\bp^2)$ of rank 1 conics and contracting the
exceptional divisor down to the locus $\Omega^2 = S^2 \bp^2$ of rank 2
(dual) conics.

The image of $\Delta$, and hence of $\pi_V|_W$, is therefore the
complement of $\Omega^2$ and we are done.  \qed

\section{Genus 4}

To any nonhyperelliptic curve of genus 4 one can associate in a
canonical way a nodal cubic threefold $\tt \subset \bp^4$ as follows
(see \cite{D}).  The canonical curve $C\subset \bp^3$ lies on a unique
quadric surface $Q\subset \bp^3$ and is base-locus of a 4-dimensional
linear system of cubics; we define $\tt \subset \bp^4$ to be the image
of the rational map $$ \lambda_{|\iic(3)|} : \bp^3 \rightarrow \bp^4.
$$ Let us denote by $\trig,\trigg \in \Theta \subset \pic^3(C)$ the two
trigonal line bundles on the curve. These satisfy $\trig\otimes \trigg
= K$ and coincide precisely when the curve has a vanishing theta-null.
The quadric surface $Q$ is ruled by trisecants $\overline{D}\subset
\bp^3$ of the curve, where $D$ belongs to $|\trig|$ or $|\trigg|$ (and
$Q$ is singular precisely when the two pencils coincide); it therefore
contracts to a point $t_0 \in \tt$ under $\lambda_{|\iic(3)|}$. Any
hyperplane through $t_0 \in \bp^4$ then pulls back to $Q$ plus a
residual hyperplane, and it follows that projection away from the point
$t_0$ is the birational inverse of $ \lambda_{|\iic(3)|}$:

\begin{equation} \begin{array}{crccl} &&\tt&\subset & \bp^4\\ &&&&\\
&\hidewidth\lambda_{|\iic(3)|}\nearrow&&\searrow&\downarrow \pi_{t_0}\\
&&&&\\ \bp^3&&=&&\bp^3\\ \end{array} \end{equation} This allows us to
see that $\tt$ is a cubic: a general hyperplane $H\subset \bp^4 =
|\iic(3)|\dual$ identifies with $\bp^3$ under the projection
$\pi_{t_0}$, and under this identification its intersection with $\tt$
is the cubic surface corresponding to the point of $|\iic(3)|$
annihilated by $H$.

\begin{prop} $\mult_{t_0}  \tt=2 $ and the projectivised tangent cone
at this point is $\bp T_{t_0} \tt = Q \subset \bp^3$.  \end{prop}

\pf That $\mult_{t_0}  \tt=2 $ follows at once from the fact that the
projection $\pi_{t_0}: \tt \rightarrow \bp^3$ is birational and $\deg
\tt =3$. On the other hand, let $H\subset \bp^4$ be any hyperplane
passing through $t_0$ and $H' = \pi_{t_0}(H) \subset \bp^3$ its
projection. Then $H\cap \tt$ is the cubic surface obtained by blowing
up the six points $C\cap H' \subset \bp^2$; these six points lie on the
conic $Q' = Q\cap H'$ and it is well-known that the resulting cubic
surface is nodal with projectivised tangent cone $Q'\subset \bp^2$ at
the node.  \qed

\begin{rem}\rm
From this we can easily write down an equation for $\tt$: choosing a
simplex of reference with $t_0 \in \bp^4$ as one vertex, and the
opposite face corresponding to a choice of cubic surface $F\in
|\iic(3)|$, the threefold $\tt$
has equation $$ z_0 Q(z_1,\ldots ,z_4) + F(z_1,\ldots ,z_4) = 0.  $$
This is the description given by Donagi.  \end{rem}

We shall need next the Fano surface $F(\tt)$ of lines on $\tt$, which
is easy to describe. First note that there is an inclusion $$
\begin{array}{rcl} i: C &\hookrightarrow& F(\tt) \\ p&\mapsto& l_p =
\hbox{line joining $t_0$ to $p\in C\subset \bp^3$.}\\ \end{array} $$ In
other words, $l_p\subset \bp^4$ is the line joining $t_0$ to the point
$p$ on the canonical curve via the projection $\pi_{t_0}$, and it is
easy to see that these are precisely the lines through $t_0 \in \bp^4$
which lie on $\tt$.  We now map $$ \begin{array}{rcl} f: S^2 C
&\rightarrow& F(\tt) \\ p+q&\mapsto& l_{pq} = \hbox{residual line in
$\tt\cap {\rm Span}\{l_p,l_q\}$.}\\ \end{array} $$ (Note that this
still makes sense on the diagonal of $S^2 C$: if $p=q$ then ${\rm
Span}\{l_p,l_q\}$ is interpreted to mean the 2-plane spanned by $t_0$
and the tangent line to the canonical curve at $p\in C$.)

If the secant line $\overline{pq}\subset \bp^3$ is not on $Q$ then
$\lambda_{|\iic(3)|}(\overline{pq}) = l_{pq}$; whilst if $\overline{pq}
\subset Q$ then it contracts down to $t_0$, but $l_{pq} = l_r$ where
$r\in \overline{pq}\cap C$ is the third point of the trisecant. Thus
$f$ is a birational morphism and is injective away from the two curves
$C\hookrightarrow S^2 C$ defined by $r\mapsto \trig(-r)$ and $r\mapsto
\trigg(-r)$, each of which it identifies with $i(C)$:

\begin{equation} \label{pic} figure 
\end{equation}

Izadi \cite{I} makes use of the lines on $\tt$ to identify $\tt \subset
\bp^4 \idfy \g00$ (theorem \ref{izadi} below).  Namely, for $r\in C$
and for $p+q\in F(\tt)\backslash i(C)$ (which we identify with the
corresponding subset of $S^2 C$ as above) she constructs  pencils which
we shall denote by $l'_r,l'_{pq} \subset \g00$ respectively. These are
characterised by their base locus: for any $p,q \in C$ let

\begin{equation} \Sigma_{pq} = C-C \cup W_2 -p-q \cup p+q-W_2 \subset
J_C.  \end{equation} Then the pencil $l'_r\subset \g00$ has base locus
$\Sigma_{pq}\cup \Sigma_{p'q'}$ where $f^{-1}(i(r)) = \{p+q,p'+q'\}$,
i.e. $p+q+r\in |\trig|$ and $p'+q'+r \in |\trigg|$; and the pencil
$l'_{pq} \subset \g00$ has base locus $\Sigma_{pq} \cup \Sigma(X)$
where (in Izadi's notation---see \cite{I}, \S7)

\begin{equation} \label{Sigma(X)} \begin{array}{rccl} \Sigma(X) &=&
\{s+t-s'-t'\ |& s,t,s',t'\in C,\\ &&&  h^0(K-p-q-s-t)>0,\\
&&&h^0(K-p-q-s'-t')>0\}.\\ \end{array} \end{equation} (In this notation
$X$ denotes a curve of genus 5 in the fibre of the Prym map over $J_C$;
though this will not concern us here.)

Izadi's result, in our (nonhyperelliptic Jacobian) situation, can then
be stated as follows.

\begin{theo} \label{izadi} Let $C\subset \bp^3$ be a canonical curve of
genus 4, and $\bp^4 = |\iic(3)|\dual$ be the ambient space of its
associated cubic threefold $\tt$. Then there is a natural
identification $\bp^4 \idfy \g00$ under which $l_{r} \idfy l'_{r}$,
$l_{pq} \idfy l'_{pq}$ and the node $t_0\in \tt$ maps to $\Theta-\trig
\cup \Theta-\trigg$.  \end{theo}

We now return to consider the Brill-Noether loci in $\su(2,K)\subset
\bp^{15}$ and to state our main result. Recall that $\ww^2 =  \g00 \cap
\su(2,K)$ (by theorem~\ref{w2g00}).

\begin{theo} \label{w2cubic} If $C$ is a nonhyperelliptic curve of
genus 4 then $\ww^2\subset \su(2,K)\subset \bp^{15}$ is precisely the
Donagi-Izadi cubic threefold $\tt \subset \g00 = \bp^4$; with node at
$\trig \oplus \trigg \in \kum(J_C^3)$.  \end{theo}

\begin{rem}\rm \label{w3node} Note that (up to S-equivalence) $\trig
\oplus \trigg$ is the unique semistable bundle in $\su(2,K)$ with $h^0
= 4$, and so $\ww^3$ is by definition empty. This is a consequence of
Mukai's lemma~\ref{muk}: since $|\trig|$ is base-point-free, $h^0(E)
\geq 4$ would imply that $h^0({\trig}^{-1}\otimes E) \geq 1$, and hence
by semistability that $\trig \subset E$. So $E$ is S-equivalent to
$\trig \oplus K{\trig}^{-1} = \trig \oplus \trigg$.

(In fact, if $\trig \not= \trigg$ then this is an isomorphism since
$\trigg \subset E$ by the same argument. If, on the other hand, $\trig
= \trigg$ then one can check using the arguments of section 3 that as
well as $\trig \oplus \trig$ there is, up to isomorphism, a unique
nonsplit extension $E$ with $h^0(E) =4$: the space of all such
extensions is $\bp (\trig) = \bp H^0(K)\dual$, in which the canonical
curve lies on a quadric cone. $E$ is then the extension corresponding
to the vertex of the cone.)

Thus it makes sense to view $\ww^3 = \{\trig \oplus \trigg\}$.
\end{rem}

\begin{lemm} Suppose that $C$ is nonhyperelliptic of genus 4 or
nontrigonal of genus 5. Then for every stable $E\in \su(2,K)$ with
$h^0(E)=3$ the exterior multiplication map $\phi_E: \bigwedge^2 H^0(E)
\rightarrow H^0(K)$ is injective.  \end{lemm}

\pf Since every element of $\bigwedge^2 H^0(E)$ is decomposable, i.e.
of the form $s\wedge t$, a nontrivial element of $\ker \phi_E$ would
give two independent sections $s,t$ generating a line subbundle
$L\subset E$. Then $r(L) \geq 1$ while by stability $\deg L \leq g-2$,
contrary to the hypotheses on $C$.  \qed

In the genus 4 case the lemma determines a rational map (defined away
from the single point $\trig \oplus \trigg$) $ \pi: \ww^2 \rightarrow
\bp^3 = |K|\dual $ sending $E\mapsto \im \phi_E \subset |K|$.  In
proving theorem \ref{w2cubic} we shall in fact prove slightly more,
namely that the following diagram commutes (and we shall also extend
this diagram in theorem \ref{quartics} below):

\begin{equation} \label{nodalcubic} \begin{array}{rcl} \tt & = & \ww^2
\subset \g00\\ & \hidewidth\pi_{t_0} \searrow  \qquad\swarrow \pi
\hidewidth&\\ & \bp^3 &\\ \end{array} \end{equation}

\begin{prop} \label{6.8} For each $p\in C \subset \bp^3$ the closure of
the fibre $\pi^{-1}(p)\subset \ww^2$ is a Hecke line $l_F$ with $\det F
=K(p)$. Moreover, $l_F$ is the unique such Hecke line contained
in~$\ww^2$ and passing through the point $\trig \oplus \trigg$.
\end{prop}

\pf Consider a stable bundle $E\in \pi^{-1}(p)$. By definition the
sections of $E$ fail to generate $E$ at the point $p$, and we denote by
$D_p\subset E_p$ the line in the fibre at $p$ which is generated by
global sections. This line uniquely determines an extension $$
\ses{E}{F}{\c_p}, $$ and by construction the coboundary map $\delta :
\c \rightarrow H^1(E)$ vanishes, so that $h^0(F) = 4$ and
$H^0(E)\subset H^0(F)$ coincides with the subspace $H^0(F(-p))$.
Finally, since $E$ is stable $F$ is stable. This proves the first part
of the proposition, with $l_F \subset \ww^2$.

Now choose a section $s\in H^0(F)$ not lying in $H^0(E)$. Then $s(p)
\not= 0$ and spans a line in the fibre $F_p$; we consider a nonzero
homomorphism $u: F\rightarrow \c_p$ such that this lines coincides with
$\ker u_p$. Then by construction $\ker u \subset F$ is a semistable
bundle with $h^0(\ker u) =4$ and hence by remark~\ref{w3node} $\ker u =
\trig \oplus \trigg$ up to S-equivalence; and this point therefore lies
on $l_F$.

It remains to show that a Hecke line with these properties is unique.
If $\trig \not= \trigg$ then by remark \ref{w3node} $F$ has subsheaf
$\trig \oplus \trigg$, and by \cite{B1} lemma 3.2 this determines the
Hecke line $l_F$ uniquely. Alternatively one can argue similarly to the
vanishing theta-null case, to which we shall now restrict.

So assume that $\trig =\trigg$.  By theorems \ref{10} and \ref{20}
$l_F$ lies in some 4-plane $\bp(x)$, meeting the curve at the image of
a point $q\in C$: thus $x = \trig(-q)$. By the proof of \ref{10}, {\it
either} $p=q$ {\it or} $l_F$ is a trisecant $\overline D$ where $\oo(D)
=Kx^{-2}(p)$. On the other hand, the second case does not occur for the
following reason:  by (\ref{verdiagram}) and (\ref{conethm}) in section
3, $h^0 \geq 3$ in $\bp (x)$ precisely along a line $\Omega^0_x$,
projection away from which maps $C$ onto a plane conic via the linear
series $|\trig(q)|$ with the single base point $q$. Thus (since $h^0
\geq 3$ along $l_F$) $l_F = \Omega^0_x$ and meets the curve only at one
point (with multiplicity 2).  Thus $l_F = \Omega^0_{\trig-p}$ and is
uniquely determined.  \qed

\begin{rem}\rm In the case $\trig \not= \trigg$ one can show that $l_F$
is the intersection of the two 4-planes $\bp(\trig(-p))$ and
$\bp(\trigg(-p))$ (in the notation of section 1), and in each space is
the tangent line to the curve at the image of $p\in C$.  In the case
$\trig = \trigg$ just discussed in the above proof, the curve $C\subset
\bp (\trig(-p))$ has a cusp at $p\in C$. The line $l_F$, passing
through $p$, is not the tangent line but is the vertex of a rank 3
quadric containing the curve.  \end{rem}

\begin{prop} \label{6.9} For each $p\in C$ the Hecke line of the
previous proposition coincides with Izadi's pencil $l'_p = l_p =
\pi_{t_0}^{-1}(p)$.  \end{prop}

\pf Consider a stable bundle $E\in \pi^{-1}(p) = l_F$. We shall show
that $E\in l'_p$; since both sets are lines the result will follow.

So we have to show that the divisor $D_E = \{ L\in J_C | h^0(C,L\otimes
E)>0 \}$ contains the surfaces $\Sigma_{st}$ and $\Sigma_{s't'}$ where
$p+s+t \in |\trig|$ and $p+s'+t' \in |\trigg|$.  (Note that for $E\in
\ww^2$ we have $D_E\in \g00$ and hence $C-C \subset D_E$ a priori---see
section 4.) Since $D_E$ is symmetric it is enough to prove that $W_2
-s-t \subset D_E$, i.e. that $h^0(E(p+q-s-t)) >0$ for all $p,q\in C$.
We shall show that $h^0(E(-s-t))\geq 1$ (and note that by proposition
\ref{4.3} this is actually equivalent); and similarly that
$h^0(E(-s'-t'))\geq 1$.

By hypothesis $\im \phi_E = H^0(K(-p))$; and we have a natural
2-dimensional subspace $V\subset H^0(K(-p))$, namely $$ V =
H^0(K(-p-s)) = H^0(K(-p-t)) = H^0(K(-p-s-t)).  $$ So consider the
subspace $\phi_E^{-1}(V) \subset \bigwedge^2 H^0(E)$ and choose
sections $u,v,w \in H^0(E)$ such that $u\wedge v, u\wedge w$ form a
basis of $\phi_E^{-1}(V)$. Since $v\wedge w \not\in \phi_E^{-1}(V)$ the
effective divisor $(u\wedge w)\in |K|$ is not supported at $s$ or $t$;
this implies that the sections $v,w$ generate $E$ at the points $s,t\in
C$. However, by construction $s+t \leq (u\wedge v)$ and $s+t \leq
(u\wedge w)$; and we claim that this can only occur if $u(s) = u(t) =
0$. For if $u(s) \not= 0$, for example, then $\c u(s) = \c v(s)$ and
$\c u(s) = \c w(s)$ and hence $s \in {\rm supp}(v\wedge w)$, a
contradiction.  Hence we obtain a nonzero section $u\in H^0(E(-s-t))$;
and similarly we can do the same for $H^0(E(-s'-t'))$.  \qed

Let us return to the proof of theorem \ref{w2cubic}. In seeking stable
bundles with three sections we may consider extensions $E\in \bp
\ext^1(K-D,D)$ with $\deg D = 1$ or 2 (using remark \ref{nag}). If
$E\in \Omega^n \backslash \Omega^{n-1}$ then by (\ref{conethm}) $$
h^0(E) = 5 - \cliff(D) - n.  $$ Thus either $D=p\in C$, and $E\in
\Omega^1$; or $D=p+q\in S^2 C$, and $E\in \Omega^0_{p+q}\cong \bp^1$.
The second case is that of example \ref{d=g-2}; we shall show next that
this case exhausts all such bundles.

\begin{prop} $\ww^2 = \bigcup_{p+q\in S^2 C} \Omega^0_{p+q}$.
Moreover, $\Omega^0_{p+q}$ maps under $\pi: \ww^2 \rightarrow \bp^3$
onto the secant line $\overline{pq}$ if $p+q \not\in i(C)$, while
$\Omega^0_{p+q} = \pi^{-1}(r)$ if $p+q+r \in |\trig|$ or $|\trigg|$.
\end{prop}

\pf We first observe (by considering diagram (\ref{verdiagram})) that
the line $\Omega^0_{p+q} \subset \bp(p+q)$ meets the image of the curve
if and only if $f(p+q) \in i(C)$ (see (\ref{pic})); and in this case
meets the curve at a point $r\in C$ representing the bundle $\trig
\oplus \trigg$. By theorem \ref{10} $\Omega^0_{p+q}$ is a Hecke line
$l_F$, where one easily checks that $\det F = K(r)$. So by the
uniqueness statement in \ref{6.8} $\Omega^0_{p+q} = \pi^{-1}(r)$.

We may now assume, then, that $E\in\ww^2$ is a stable bundle for which
$\pi(E)$ does not lie on the canonical curve; $\pi(E)$ then lies on
some secant line $\overline{pq}\subset \bp^3$. This means that $\im
\phi_E\subset H^0(K)$ is a hyperplane, distinct from $H^0(K(-p))$ and
$H^0(K(-q))$ but containing the 2-dimensional subspace $H^0(K(-p-q))$.
As in the proof of the previous proposition we can find a basis $u,v,w
\in H^0(E)$ such that $u\wedge v, u\wedge w$ are a basis of
$\phi_E^{-1} H^0(K(-p-q))$ which $v\wedge w$ completes to a basis of
$\bigwedge^2 H^0(E)$. Then we have $$ p+q\leq (u\wedge v), \quad
p+q\leq  (u\wedge w), $$ while $p+q \not\leq (v\wedge w)$. As before,
it follows from this that $u(p) = u(q) = 0$, i.e. $h^0(E(-p-q))> 0$ and
so $E\in \Omega^0_{p+q}$.  \qed

{\it Proof of theorem \ref{w2cubic}.} By propositions \ref{6.8},
\ref{6.9} and theorem \ref{izadi} it suffices to check that the line
$\Omega^0_{p+q}\subset \ww^2 \subset \g00$ coincides with the pencil
$l'_{pq}$ if $p+q \in F(\tt)\backslash i(C)$, i.e. when
$h^0(\trig(-p-q)) = h^0(\trigg(-p-q)) = 0$.

Consider a stable bundle $E\in \Omega^0_{p+q}$. Since $h^0(E(-p-q))> 0$
the symmetric divisor $D_E$ trivially contains the surfaces $W_2 - p-q$
and $p+q - W_2$; while $C-C \subset D_E$ since $E\in \ww^2$.  We will
show that $D_E$ also contains the surface $\Sigma(X)$ (see
(\ref{Sigma(X)})). Let $\lambda = \oo(s+t-s'-t') \in \Sigma(X)$. By
definition we have an exact sequence $$ \ses{\lambda(p+q)}{\lambda
\otimes E}{K\lambda (-p-q)} $$ with, say, extension class $f\in
\ext^1(K-p-q,p+q) = H^0(C,K^2(-2p-2q))\dual$.

We have to show that $h^0(\lambda \otimes E) >0$. We can suppose that
$h^0(\lambda(p+q)) = 0$ (otherwise there is nothing to prove); so by
Riemann-Roch $h^1(\lambda(p+q)) = 1$. If $h^0( K\lambda (-p-q)) >1$
then $h^0(\lambda \otimes E) >0$ and we are done; so we assume $h^0(
K\lambda (-p-q)) =1$. In this case $h^0(\lambda \otimes E) >0$ if and
only if the coboundary map $$ \delta: H^0(K\lambda (-p-q)) \rightarrow
H^1(\lambda(p+q)) $$ vanishes, which in turn is equivalent to $\ker f$
containing the image of the multiplication map:  $$ H^0(K\lambda
(-p-q)) \otimes H^0(K\lambda^{-1} (-p-q)) \rightarrow \ker f \subset
H^0(K^2 (-2p-2q)).  $$ In fact we shall check that the image is
contained in the subspace $S^2 H^0(K(-p-q)) \subset \ker f$ (see
example \ref{d=g-2}).

This last assertion results from the definition (\ref{Sigma(X)}): we
can write $$ K = \oo(p+q+s+t+u+v) = \oo(p+q+s'+t'+u'+v'), $$ for some
$u,v,u',v' \in C$, and hence $$ \begin{array}{rclcl} K\lambda(-p-q) &=&
K(-p-q-s'-t'+s+t) &=& \oo(u'+v'+s+t), \\ K\lambda^{-1}(-p-q) &=&
K(-p-q-s-t+s'+t') &=& \oo(u+v+s'+t'). \\ \end{array} $$ By hypothesis
these divisors are unique in their linear equivalence classes and we
can write their sum as $$ (u'+v'+s+t)+(u+v+s'+t') =
(s+t+u+v)+(s'+t'+u'+v') $$ where $s+t+u+v$ and $s'+t'+u'+v' \in
|K(-p-q)|$.  \qed

\medskip

We shall conclude this section by giving another interpretation of
diagram (\ref{nodalcubic}), as follows. First, we shall view
$\Gamma_{00} \hookrightarrow S^4 H^0(C,K)$ by assigning to each element
the leading terms of its Taylor expansion at $0\in J_C$; or
equivalently by assigning to a divisor its tangent cone at the origin.
This map is injective by \cite{I} lemma 2.1.1.

Next we note that there is a distinguished element $q^2 \in S^4
H^0(C,K)$, where $q\in S^2 H^0(C,K)$ is the equation of the quadric
$Q\subset \bp^3$ containing the canonical curve. Under the above
inclusion this comes from the split divisor $\Theta-\trig \cup
\Theta-\trigg \in \g00$.

Third, we identify $\bp^3 = \bp T_0 J_C$ with the space of
translation-invariant vector fields on the Jacobian. One can then map
$$ \begin{array}{rcrcl} \alpha &:& \bp T_0 J_C &\rightarrow& \bp
\bigl(S^4 H^0(K) / \c q^2 \bigr)\\ &&D &\mapsto & qDf - fDq,\\
\end{array} $$ where $f\in H^0(\iic(3))$ is any cubic through the
canonical curve. It is easy to check that this construction is
independent of the choice of $f$; moreover $\alpha$ is an isomorphism
onto the subspace $\bp (\Gamma_{00} / \c q^2)$, as observed by
Beauville--Debarre \cite{BD}, pages 32--33.

\begin{theo} \label{quartics} The following diagram commutes:  $$
\begin{array}{ccccc} \ww^2 & \subset &  \g00 &\subset & \bp S^4
H^0(K)\\ &&&&\\ \pi \downarrow &&\downarrow&&\downarrow\\ &&&&\\ \bp^3
& \map{\alpha}& \bp \bigl(\Gamma_{00} / \c q^2\bigr)&\subset&
\bp\bigl(S^4 H^0(K) / \c q^2 \bigr)\\ \end{array} $$ \end{theo}

\pf We have to check commutativity of the left-hand square, and since
both vertical arrows are linear projections it is sufficient to check
commutativity over points $p\in C$ of the canonical curve. For such a
point denote by $D_p\in \bp T_0 J_C$ the associated constant vector
field. By propositions \ref{6.8} and \ref{6.9} the line $\pi^{-1}(p)$
corresponds to the pencil $l_p'$ with base locus $\Sigma_{st}\cup
\Sigma_{s't'}$, where the points $s,t,s',t' \in C$ are defined by
$p+s+t \in |\trig|$ and $p+s'+t' \in |\trigg|$. By tangent cones at the
origin, the pencil $l'_p$ corresponds to a pencil of quartics spanned
by $qDf -fdq$ and $q^2$, for some $D\in T_0 J_C$ and $f\in
H^0(\iic(3))$. We have to show that $D=D_p$.

Since the pencil is uniquely determined by (the tangent cone of) its
base locus, it is enough to check that the two quartics $qD_pf -fD_pq$
and $q^2$ contain the tangent cones of $C-C$, $W_2 -s-t$ and $W_2
-s'-t'$. These tangent cones are the canonical curve $C\subset \bp^3$
and the two trisecants $\overline{pst}, \overline{ps't'} \subset \bp^3$
respectively. The result now follows easily: $q$ vanishes on all three
curves; while $f$ (and hence $fD_p q$) vanishes on $C$, and---since the
two trisecants span the tangent plane to $Q$ at $p$---the derivative
$D_p q$ vanishes on the two lines.  \qed

\section{Genus 5}

Let $C$ be a curve of genus 5. If $C$ is nontrigonal then the canonical
curve $C\hookrightarrow \bp^4$ is the complete intersection of a net of
quadrics $|\iic(2)| = \bp^2$, in which the locus $\Gamma \subset \bp^2$
of singular quadrics is a plane quintic curve, smooth if $C$ has no
vanishing theta-nulls, otherwise having ordinary double points
corresponding to quadrics of rank 3 (see \cite{ACGH}, page 270).

Let $\thsing = W^1_4$ be the singular locus of the theta divisor. This
is a curve, and by assigning to each point $x\in \thsing$ its
projectivised tangent cone $\bp T_x \Theta = Q_x$ we have a double
cover $$ \begin{array}{rrcl} f:&\thsing & \rightarrow & \Gamma \subset
\bp^2 \\
   &  x & \mapsto & Q_x \\ &&& \displaystyle = \bigcup_{D\in |x|}
\overline{D} \subset \bp^4.\\ \end{array} $$ The sheet interchange of
$\thsing$ with respect to this double cover is induced by the Serre
involution of $J_C^4$.

\begin{lemm} $f^* \oo_{\Gamma}(1) = \oo_{\thsing} (\Theta)$. Moreover,
the induced restriction map $H^0(J_C^{g-1}, 2\Theta) \rightarrow
H^0(\Gamma, \oo(2))$ is surjective.  \end{lemm}

\pf The first part follows from \cite{G}.  To prove that the pull-back
of hyperplane sections is surjective, it is sufficient to show this on
the image of the Kummer map $J_C \rightarrow |2\Theta|$, $a\mapsto
\Theta_a + \Theta_{-a}$. In other words, we consider the rational map
$$ \alpha : J_C \rightarrow |\oo_{\Gamma}(2)| \cong \bp^5 $$ sending
$a\in J_C$ to the divisor whose pull-back to $\thsing$ is $(\Theta_a +
\Theta_{-a})\cap \thsing$. (Note that $\alpha$ is defined away from
$C-C \subset J_C$: this follows from \cite{W}, theorem 2.4.) One can
show that the map $$ \begin{array}{rcl} \beta:  J_C &\rightarrow &
S^{10}(\thsing) \\ a &\mapsto& \Theta_{a}\cap \thsing\\ \end{array} $$
is injective (see, for example, \cite{ACGH}, pages 265--268); this
implies that $\alpha$ is a finite map and so we are done.  \qed

It follows from this that the image of $\thsing$ under the Kummer map
is a Veronese embedding of $\Gamma \subset \bp^2$:

\begin{equation} \label{mapv} \begin{array}{ccc} \thsing & \map{\kum} &
\bp^{31} \\ f\downarrow && \uparrow v \\ \Gamma & \subset & \bp^2 \\
\end{array} \end{equation} where $v(\bp^2)\subset \bp^5 \subset
\bp^{31}$ is a Veronese surface.

\begin{theo} \label{g5ver} For any curve $C$ of genus 5 the
Brill-Noether locus $\ww^3\subset \su(2,K)$ is a Veronese surface
intersecting the Kummer variety in the Veronese image of a plane
quintic $\Gamma \subset \bp^2$. In particular:  \begin{enumerate} \item
If $C$ is nontrigonal then $\ww^3 = v(\bp^2)$, where $\Gamma$ is as in
(\ref{mapv}) and $v(\Gamma) = \kum(\thsing)$.  \item If $C$ is trigonal
then $\Gamma \subset  \bp H^1(g^1_3)$, where $g^1_3$ is the (unique)
trigonal line bundle, is the projection of the canonical curve away
from a trisecant; and its Veronese image cuts the Kummer in the
component $C+g^1_3$ of~$\thsing$.  \end{enumerate} \end{theo}

{\it Proof of part 2.} This is easily dispatched. We first remark that
it is well-known that on a curve of genus $\geq 5$ a $g^1_3$ is unique
if it exists; while for a curve of genus 5 the following argument will
give another proof of this fact.

Let $|D| = g^1_3$; by lemma \ref{muk} any stable bundle  $E\in
\su(2,K)$ with $h^0(E) \geq 4$ has line subbundle $\oo(D) \subset E$,
so $E$ belongs to the 5-plane $\bp(g^1_3)$ of section 1. By example
\ref{d=g-2} we have seen that $h^0(E) = 4$ precisely along a Veronese
surface in $\bp(g^1_3)$.  This intersects the Kummer precisely in the
image of the curve---that is, in the Kummer image of $C+g^1_3$---and
from diagram (\ref{verdiagram}) this is the projection of the canonical
curve as asserted.  \qed

From now on we shall assume that the curve $C$ is nontrigonal.
Before proving part 1 of the theorem we shall need to make some further
observations about the curve $\Gamma$; we consider the map

\begin{equation} \label{mapl} \begin{array}{rcl} l: S^2 C &\rightarrow&
(\bp^2)\dual \\ D & \mapsto & |{\cal I}_{C\cup \overline{D}}(2)|.\\
\end{array} \end{equation} In other words $l(D)$ is the pencil of
quadrics containing $C$ and the line $\overline{D}$. Note that the base
locus of such a pencil is a quartic del Pezzo surface containing
sixteen lines, and so $\deg l = 16$.

For each $D\in S^2 C$ we shall identify the five quadrics

\begin{equation} \label{lDmeetG} l(D) \cap \Gamma = \{ Q_1, \ldots ,
Q_5\}.  \end{equation} Projection away from the line $\overline{D}
\subset \bp^4$ maps the canonical curve $C$ to a 5-nodal plane sextic
$C'\subset \bp^2$. (Note, again, that the del Pezzo base locus of the
pencil $l(D)$ is obtained by blowing up $\bp^2$ in the five nodes of
$C'$.) Let us denote by $D^{(1)}, \ldots , D^{(5)} \in S^2 C$ the
divisors over the five nodes of $C'\subset \bp^2$.  Then by
Riemann-Roch each $|D+D^{(i)}|$ is a $g^1_4$, and hence each

\begin{equation} \label{lDmeetG'} Q_i = Q_{D+D^{(i)}} =\bigcup_{D' \in
|D+D^{(i)}|}\overline{D'} \qquad i=1,\ldots , 5, \end{equation} is a
quadric of rank $\leq 4$ containing the line $\overline{D}$. These are
therefore the points of intersection (\ref{lDmeetG}).

We now return to the proof of theorem \ref{g5ver}. We consider stable
extensions $E\in \bp \ext^1(K-D,D)$ where (by remark \ref{nag} (i)) we
may take $\deg D = 2$ or 3. For such an extension, by (\ref{conethm}),

$$ h^0(E) = 6-\cliff(D) -n $$ where $E\in \Omega_D^n\backslash
\Omega^{n-1}_D$. So for $h^0(E) =4$ we must have $n+ \cliff(D) =2$; if
$\deg D =3$ then this forces $|D| = g^1_3$ contrary to the hypothesis
that $C$ is nontrigonal. So the only possibilities we need to consider
are $D\in S^2 C$, and then $h^0(E) =4$ for $E\in \Omega_D^0 \subset \bp
\ext^1(K-D,D)$.

In this situation diagram (\ref{verdiagram}) becomes:

$$ \begin{array}{rcr} C & \map{|K-D|} & C'\subset \bp^2 \\ &&\\
\scriptstyle|2K-2D| \downarrow &&\ver \downarrow  \\ &&\\ \Omega_D^0 =
\bp^1 \subset \bp^7 & \map{\delta} & \bp^5 \\ \end{array} $$ where
$C'\subset \bp^2$ is the 5-nodal sextic as above; and in particular
$\delta$ is surjective.

It follows that $\ww^3 = \bigcup_{D\in S^2 C} \eee_D\Omega_D^0 \subset
\bp^{31}$, where we observe that {\it each $\eee_D \Omega_D^0$ is a
nonsingular conic.} This is because by lemma \ref{bertram} the rational
map $\eee_D$ comes from the (complete) linear series $|\iic (2)|$ on
$\bp^7$; while $\Omega_D^0$ has no intersection with the base locus $C$
since $|K-D|$ has no base points---i.e. the canonical curve has no
trisecant lines since $C$ is nontrigonal.

Finally, theorem \ref{g5ver} will follow directly once we observe that

\begin{equation} \label{Dconic} \eee_D \Omega_D^0 = v(l(D)) \qquad
\hbox{for all $D\in S^2 C$,} \end{equation} where $v$ and $l$ are as
defined in (\ref{mapv}) and (\ref{mapl}) respectively.  To prove
(\ref{Dconic}) it is sufficient, since both sides are nonsingular
conics, to show that both contain the five points $v(Q_1),\ldots ,
v(Q_5)$ (see (\ref{lDmeetG}) and (\ref{lDmeetG'})); and so it remains
only to check this for $\eee_D \Omega_D^0$.

First note that an extension $E\in \Omega_D^0$ fails to be stable (i.e.
maps to the Kummer variety) if and only if it lies on $\secant^2 C$;
and there are precisely five such points, which are the intersections
of $\Omega_D^0$ with the secant lines
$\overline{D^{(1)}},\ldots,\overline {D^{(5)}}$, where as before the
$D^{(i)}\in S^2 C$ are the nodal divisors over the curve $C'\subset
\bp^2$. The corresponding extensions then contain $\oo(K-D-D^{(i)})$,
respectively, as line subbundles---in other words they map under
$\eee_D$ to the points

$$ \oo(D+D^{(i)}) \oplus \oo(K-D-D^{(i)}) \in \kum(J_C^{g-1}).  $$ By
(\ref{lDmeetG'}) and diagram (\ref{mapv}) these are precisely the
images $v(Q_1),\ldots,v(Q_5)$, which completes the proof.  \qed

\section{Genus 6}

By proposition \ref{bound} we have $h^0(E) \leq 6$ for all semistable
bundles $E$ in $\su(2,K)$ on a nonhyperelliptic curve $C$ of genus 6,
and $h^0(E) \leq 5$ if $C$ is not trigonal or a plane quintic.

\begin{theo} Let $C$ be a nonhyperelliptic curve of genus 6.
\begin{enumerate} \item If $C$ is not trigonal or a plane quintic then
there exists a unique stable bundle $E\in \su(2,K)$ with $h^0(E) = 5$;
i.e. $\ww^4 = \{E\}$.  \item If $C$ is trigonal then $\ww^4 \cong
\bp^1$ is a line, along which $h^0(E) = 5$, and does not meet
$\kum(J_C^{g-1})$.  \item If $C$ is a plane quintic then $h^0(E)\geq 5$
if and only if $E$ is in the S-equivalence class of the point $g^2_5
\oplus g^2_5 \in \kum(J_C^{g-1})$.  \end{enumerate} \end{theo}

\pf First suppose that $C$ is not a plane quintic, and observe that if
$h^0(E)\geq 5$ and $E$ is semistable then it is necessarily stable:
otherwise $E$ fits in an extension (\ref{basicext})  with $\deg D = 5$
(and with $D$ not necessarily effective). Then $h^0(D) = h^0(K-D) \leq
2$, since $C$ is not a plane quintic, and so $h^0(E) \leq 4$.

On the other hand, for any tetragonal pencil $g^1_4 = |D|$ one can
(following Mukai) apply lemma \ref{muk} to observe that $h^0(E(-D))
\geq 1$ for any such bundle. By stability this means that $\oo(D)
\subset E$ is a line subbundle, i.e. $E \in \bp (D)$, the corresponding
6-plane of the ruling of section 1. By example \ref{d=g-2}, $h^0(E) =5$
exactly for $E\in \Omega^0$; if $C$ is nontrigonal this is a single
point, and part 1 is proved.

If $C$ is trigonal then we may take $D=K-2L$ where $|L| = g_3^1$;
$|K-D|$ maps $C\rightarrow \bp^2$ with degree 3 onto a conic, so in
this case $\ker \delta$ is 2-dimensional and $\Omega_D^0 \subset
\bp(D)$ is a line. This line does not meet the image of $C$ in $\bp(D)$
since $|K-D|$ is base-point-free; so we have proved part 2.

For part 3, first note that by the reasoning of remark \ref{w3node} the
only semistable bundle with $h^0(E) =6$ is $E=g^2_5 \oplus g^2_5$. On
the other hand the reasoning of part 1 above yields extensions with
$h^0(E) =5$ in $\Omega^0_D$ for any $|D| = g^1_4$. In this case the
tetragonal pencils are precisely the projections from points of the
plane quintic, i.e. $D=L-p$ for $|L| = g^2_5$ and some $p\in C$. The
map $C\rightarrow \bp^2$ given by the series $|K-D|$ is projection of
the canonical curve (which lies on a Veronese surface) away from the
conic in $\bp^5$ spanned by $D$, and hence has base-point $p$. The
image is thus the plane quintic model of $C$; in particular $\delta$ is
surjective and $\Omega^0_D$ is a single point. But because $p\in C$ is
a base-point of $|K-D|$, the curve passes through $\Omega^0_D$ at the
image of $p$, which is the equivalence class of $g^2_5 \oplus g^2_5$.
This proves part 3.  \qed

\begin{rem}\rm Part 1 of the theorem was observed by Mukai in
\cite{Muk}. Recall that the canonical curve lies on the (del Pezzo)
transverse intersection with a $\bp^5 \subset \bp^9$ of the Pl\"ucker
embedded Grassmannian of lines in $\bp^4$. The bundle $E$ in the
theorem is then dual to the restriction to the curve of the
tautological bundle on the Grassmannian.  \end{rem}

It is well-known that a generic curve of genus 6 possesses five
tetragonal pencils, so the point $E$ of part 1 is common to the
corresponding five 6-planes of the ruling. It is amusing to see this
using the results of section~1.

Let us denote the five by $x_0,\ldots , x_4 \in \pic^4(C)$; and let us
recall how they are related to each other. By Riemann-Roch $|x_0| =
g^1_4$ implies that $|Kx_0^{-1}| = g^2_6$. Thus the image of $$
\lambda_{|Kx_0^{-1}|} : C \rightarrow \bp^2 $$ is a sextic with four
nodes, which we shall denote by $p_1,\ldots,p_4 \in \bp^2$:

\begin{center} 
figure 
\end{center} 
Let $D_i\in S^2C$,
$i=1,\ldots,4$, be the nodal divisors, i.e.  $p_i =
\lambda_{|Kx_0^{-1}|}(D_i)$. If we denote by $H= Kx_0^{-1}$ the
hyperplane class on $C$ then by adjunction in the blow-up at the four
nodes we have $K = 3H - D_1-\cdots -D_4$ and hence $$ x_0 = \oo(2H -
D_1 - \cdots -D_4) $$ i.e. {\it $|x_0|$ is cut out by the pencil of
conics through the four nodes $p_1, \ldots ,p_4$}. In this model it is
easy to see the remaining four $g^1_4$s: {\it for $i=1,\ldots , 4$ the
pencil $|x_i|$ is cut out by the lines through $p_i \in \bp^2$}.
Formally $x_i = \oo(H-D_i)$, and in particular we deduce that

\begin{equation} \label{tets} x_0 \otimes x_i = \oo(K - D_i).
\end{equation}

Consider again the five 6-planes $\bp(x_0),\ldots , \bp(x_4)$: we have
just seen that the bundles $E$ for which $h^0(E)=5$ are the points
$\Omega^0$ in these five spaces. Let us denote these five bundles by
$E_i \in \bp(x_i)$. We are claiming that they all coincide:

\begin{equation} \label{5=} E_0 = E_1 = E_2 = E_3 = E_4 \in \su(2,K).
\end{equation} To see this, let us work in $\bp(x_0)$. In example
\ref{d=g-2} we have seen that $\Omega^0 = \{E_0\}$ is the vertex of a
Veronese cone $\Omega^1$ containing the image of $C$ (i.e. the
intersection of $\bp(x_0)\subset \su(2,K)$ with $\kum(J_C^{g-1})$), and
that projection away from the vertex maps $C$ to $\bp^2$ via the linear
system $|Kx_0^{-1}|$. The image of this map is the 4-nodal sextic just
noted, and it follows that the four secant lines $\overline{D_i}
\subset \bp(x_0)$, $i=1,\ldots ,4$, where $D_i \in S^2 C$ is the $i$-th
nodal divisor as above, all pass through the vertex~$E_0$.  By
proposition \ref{40} together with (\ref{tets}) it follows that $$
\overline{D_i} = \bp (x_0) \cap \bp(x_i), \qquad i=1,\ldots , 4.  $$
Thus $E_0 \in \overline{D_i}\subset \bp(x_i)$ for each $i$, and
therefore coincides with the unique bundle $E_i \in \bp(x_i)$ having
five sections---so again we have proved (\ref{5=}).

\begin{rem}\rm \label{21dimspan} By the proof of theorem \ref{w2g00} we
have, for a curve of genus 6, $ \ww^2 \subset \g00 \cap \su(2,K)
\subset \ww^2 \cup \bigcup_{x\in W^1_4} \bp(x).  $ We have seen that in
each such $\bp(x)$, $|x| = g^1_4$, we have $h^0(E)\geq 3$ along a cubic
cone, which in particular spans $\bp(x)$. It follows that $$ \g00 \cap
\su(2,K) = \ww^2 \cup \bigcup_{x\in W^1_4} \bp(x), $$ and that this
intersection properly contains $\ww^2$ since $h^0(E) =2$ at the generic
point of each $\bp(x)$.  For $C$ generic $\bigcup_{x\in W^1_4} \bp(x)$
consists of five 6-planes meeting pairwise in ten lines concurrent at
the point $\ww^4$.  \end{rem}

\bigskip

\noindent {\addressit Department of Mathematical Sciences, Science
Laboratories, South Road, Durham DH1~3LE, U.K.  }

\noindent {\addressit E-mail:} {\eightrm w.m.oxbury@durham.ac.uk}

\medskip \noindent {\addressit DPMMS, University of Cambridge, 16 Mill
Lane, Cambridge CB2~1SB, U.K.  }

\noindent {\addressit E-mail:} {\eightrm pauly@pmms.cam.ac.uk}

\medskip \noindent {\addressit Department of Mathematics, Boston
University, Boston, MA~02215, U.S.A.  }

\noindent {\addressit E-mail:} {\eightrm ep@math.bu.edu} \medskip

\begin{thebibliography}{ACGH}

\bibitem{ACGH}  {\capit E. Arbarello, M. Cornalba, P.A. Griffiths, J.
Harris}: {\sl Geometry of Algebraic Curves}, Springer 1985;


\bibitem{B1} {\capit A. Beauville}: Fibr\'e de rang deux sur une
courbe, fibr\'e d\'eterminant et fonctions theta II, Bull. Soc. Math.
France, 119 (1991) 259--291;


\bibitem{BD} {\capit A. Beauville, O. Debarre}: Sur les fonctions theta
du second ordre, {\sl Arithmetic of Complex Manifolds}, Lect. Notes
Math. 1399 (1989) 27--39;

\bibitem{BNR} {\capit A. Beauville, M.S. Narasimhan, S. Ramanan}:
Spectral curves and the generalised theta divisor, J. Reine Angew.
Math. 398 (1989) 169--179;

\bibitem{B} {\capit A. Bertram}: Moduli of rank 2 vector bundles, theta
divisors and the geometry of curves in projective space, J. Diff. Geom.
35 (1992) 429--469;

\bibitem{BF} {\capit A. Bertram, B. Feinberg}: On stable rank two
bundles with canonical determinant and many sections, to appear in
Proceedings of the Europroj Annual Conferences, Catania/Barcelona
1993/4;

\bibitem{C} {\capit A. Coble}: {\sl Algebraic Geometry and Theta
Functions}, AMS Colloquium Publications vol X, 1929;

\bibitem{Deb} {\capit O. Debarre}: Sur la d\'emonstration de A. Weil du
th\'eor\`eme de Torelli pour les courbes, Compositio Math. 58 (1986)
3--11;

\bibitem{D} {\capit R. Donagi}: The fibres of the Prym map, in {\sl
Curves, Jacobians and Abelian Varieties}, Contemporary Mathematics vol
136 (1992);

\bibitem{vGvdG} {\capit B. van Geemen, G. van der Geer}: Kummer
varieties and the moduli spaces of abelian varieties, Am. J. Math. 108
(1986) 615--642;

\bibitem{G} {\capit M. Green}: Quadrics of rank four in the ideal of
the canonical curve, Invent. Math. 75 (1984) 85--104;

\bibitem{GL} {\capit M. Green, R. Lazarsfeld}: On the projective
normality of complete linear series on an algebraic curve, Invent.
Math. 83 (1986) 73--90;


\bibitem{I} {\capit E. Izadi}: The geometric structure of $\aa_4$, the
structure of the Prym map, double solids and $\Gamma_{00}$-divisors, J.
Reine Angew. Math. 462 (1995) 93--158;

\bibitem{LN} {\capit H. Lange, M.S. Narasimhan}: Maximal subbundles of
rank two vector bundles on curves, Math. Ann. 266 (1983) 55--72;

\bibitem{Lan2} {\capit H. Lange}: Higher secant varieties of curves and
the theorem of Nagata on ruled surfaces, Manuscripta Math. 47 (1984)
263--269;


\bibitem{Lan} {\capit H. Lange}: H\"ohere Sekantenvariet\"aten und
Vektorb\"undel auf Kurven, Manuscripta Math. 52 (1985) 63--80;

\bibitem{L} {\capit Y. Laszlo}: Un th\'eor\`eme de Riemann pour les
diviseurs theta sur les espaces de modules de fibr\'es stables, Duke
Math. J. 64 no. 2 (1991) 333--347;

\bibitem{Muk} {\capit S. Mukai}: Curves and Grassmannians, in {\sl
Algebraic Geometry and Related Topics}, eds. J-H. Yang, Y. Namikawa, K.
Ueno, 1992;

\bibitem{M} {\capit D. Mumford}: {\sl Curves and their Jacobians}, Ann
Arbor, 1975;


\bibitem{NR} {\capit M.S. Narasimhan, S. Ramanan}: $2\theta$-linear
systems on abelian varieties, in {\sl Vector Bundles on Algebraic
Varieties}, Tata Institute, Bombay 1984;

\bibitem{PR} {\capit K. Paranjape, S. Ramanan}: On the canonical ring
of a curve, in {\sl Algebraic Geometry and Commutative Algebra Vol II},
Kinokuniya, Tokyo 1988;

\bibitem{W} {\capit G.E. Welters}: The surface $C-C$ on Jacobi
varieties and second order theta functions, Acta. Math. 157 (1986)
1--22.








\end{thebibliography}
\end{document}